\documentclass[11pt,reqno]{amsart}
\pdfoutput=1

\usepackage{latexsym,amsthm,amssymb,amsmath,amscd,hyperref}
\usepackage{tikz-cd}
\usetikzlibrary{matrix,arrows,decorations.pathmorphing}
\usepackage{bm}
\usepackage{fancyhdr}
\usepackage{mathtools}

\usepackage{tcolorbox}

\usepackage[normalem]{ulem}

\usepackage[hmargin=2cm,vmargin=1.5cm]{geometry}

\usepackage{amsmath}
\usepackage{amscd}
\usepackage{amssymb}
\usepackage{latexsym}

\usepackage{graphicx}

\newtheorem{theorem}{Theorem}[section]
\newtheorem{proposition}[theorem]{Proposition}
\newtheorem{lemma}[theorem]{Lemma}
\newtheorem{corollary}[theorem]{Corollary}
\theoremstyle{definition}
\newtheorem{definition}[theorem]{Definition}
\newtheorem{remark}[theorem]{Remark}

\newcommand{\thmref}[1]{Theorem~\ref{#1}}
\newcommand{\secref}[1]{Section~\ref{#1}}
\newcommand{\proref}[1]{Proposition~\ref{#1}}
\newcommand{\lemref}[1]{Lemma~\ref{#1}}
\newcommand{\corref}[1]{Corollary~\ref{#1}}
\newcommand{\remref}[1]{Remark~\ref{#1}}

\usepackage{comment}
\def\tcb{}

\newcommand{\Cl}{\textup{Cl}}

\newcommand{\Aut}{\textup{Aut}}
\newcommand{\Gal}{\textup{Gal}}

\newcommand{\spn}{\textup{span}}

\newcommand{\Tr}{\textup{Tr}}

 \newcommand{\CC}{\mathbb{C}}
 \newcommand{\QQ}{\mathbb{Q}}
 \newcommand{\RR}{\mathbb{R}}
 \newcommand{\NN}{\mathbb{N}}
  \newcommand{\TT}{\mathbb{T}}
 \newcommand{\ZZ}{\mathbb{Z}}
 
 \newcommand{\B}{\mathcal{B}}

\newcommand{\E}{\mathcal{E}}
\newcommand{\F}{\mathcal{F}}
\newcommand{\G}{\mathcal{G}}
\renewcommand{\H}{\mathcal{H}}
\renewcommand{\O}{\mathcal{O}}
\newcommand{\I}{\mathcal{I}}

\newcommand{\M}{\mathcal{M}}
\newcommand{\N}{\mathcal{N}}
\renewcommand{\P}{\mathcal{P}}
\newcommand{\R}{\mathcal{R}}
\newcommand{\U}{\mathcal{U}}

\renewcommand{\a}{\mathfrak{a}}
\renewcommand{\b}{\mathfrak{b}}
\renewcommand{\c}{\mathfrak{c}}
\renewcommand{\d}{\mathfrak{d}}
\renewcommand{\k}{\mathfrak{k}}
\newcommand{\m}{\mathfrak{m}}
\newcommand{\n}{\mathfrak{n}}
\newcommand{\p}{\mathfrak{p}}

\newcommand{\Ind}{\textup{Ind}}
\usepackage{xcolor}
\newcommand{\Res}{\textup{Res\;}}
\newcommand{\reg}{\textup{reg}}
\newcommand{\supp}{\textup{supp}}
\newcommand{\X}{\mathcal{X}}

\newcommand{\ex}{\mathcal{E}}
\newcommand{\res}{\textup{res}}

\newcommand{\Kz}{\mathbb{K}}

\newcommand{\Ad}{\textup{Ad}}
\newcommand{\tors}{\textup{tor}}

\newcommand{\inv}{^{-1}}

 \renewcommand\Re{\operatorname{Re}}
  \renewcommand\Im{\operatorname{Im}}
  \newcommand\sign{\operatorname{sign}}
  \def\rxrx{R\rtimes R^\times}
  \def\sofa{t_\a}
  \def\rk{\textup{rank}}
  \def\charctr{\gamma}
  \def\W{W}
  \def\unk{\underline{N}(\k)}
  \mathchardef\dash="2D
  
\begin{document}

\title[C*-dynamical invariants and actions of congruence monoids]{Partition functions as C*-dynamical invariants\\ and actions of congruence monoids}

\date{8 July 2020}

\author[Bruce]{Chris Bruce}
\address[Chris Bruce]{School of Mathematical Sciences, Queen Mary University of London, Mile End Road, E1 4NS London, United Kingdom, and
	School of Mathematics and Statistics, University of Glasgow, University Place, Glasgow G12 8QQ, United Kingdom}
\email[Bruce]{Chris.Bruce@glasgow.ac.uk}

\author[Laca]{Marcelo Laca}
  \address[Marcelo Laca]{
  Department of Mathematics and Statistics\\
  University of Victoria\\
  Victoria, BC V8W 3R4\\
  Canada}%
  \email[Laca]{laca@uvic.ca}
  
  \author[Takeishi]{Takuya Takeishi}
  \address[Takuya Takeishi]{
  Faculty of Arts and Sciences\\
  Kyoto Institute of Technology\\
  Matsugasaki, Sakyo-ku, Kyoto\\
  Japan
  }
 \email[Takeishi]{takeishi@kit.ac.jp}
\subjclass[2010]{Primary 46L55; Secondary 11M55, 11R04, 82B10.}
\thanks{Research partially supported by the Natural Sciences and Engineering Research Council of Canada through an Alexander Graham Bell CGS-D award (Bruce) and a Discovery Grant (Laca), and by JSPS KAKENHI Grant Number JP19K14551 (Takeishi).}

\begin{abstract}
We study KMS states for the C*-algebras of $ax+b$-semigroups of algebraic integers in which the multiplicative part is restricted to a congruence monoid, as in recent work of Bruce. We realize the extremal low-temperature KMS states as generalized Gibbs states in concrete representations induced from extremal traces of certain group C*-algebras. We use these representations to compute the type of extremal KMS states and we determine explicit partition functions for those of type I. The resulting collection of partition functions is an invariant for equivariant isomorphism classes of  C*-dynamical systems, which produces further invariants through the analysis of the topological structure of the KMS state space.  We use this to characterize several features of the underlying number field and congruence monoid. In most cases our systems have infinitely many type I factor KMS states and at least one type II factor KMS state at the same inverse temperature and there are infinitely many partition functions. In order to deal with this multiplicity, we establish, in the context of general C*-dynamical systems, a precise way to associate partition functions to extremal  type I KMS states. This discussion of partition functions for {C*-dynamical} systems may be of interest by itself and is likely to have applications in other contexts so we include it in a self-contained initial section that is partly expository and is independent of the number-theoretic background.
\end{abstract}

\maketitle

\setlength{\parindent}{0cm} \setlength{\parskip}{0.5cm}

\section{Introduction}
 C*-dynamical systems of number-theoretic origin have attracted sustained interest ever since Bost and Connes \cite{BC} introduced their quantum statistical mechanical system based on a noncommutative Hecke C*-algebra associated to the inclusion of rings $\ZZ \subseteq \QQ$. Their system has two remarkable features:  it exhibits spontaneous symmetry breaking  of a Galois symmetry at low-temperature and its partition function is  the Riemann zeta function. This second feature had been foreshadowed by the systems proposed independently by B. Julia \cite{Jul} and D. Spector \cite{Spe}, in which the distribution of prime numbers was interpreted in terms of the distribution of energy levels of a quantum statistical mechanical system consisting of a free bosonic gas. In this {\em Riemann gas} the non-interacting elementary particles are the prime numbers  and the energy of a prime $p$ is given by $\log p$, so that the resulting partition function is the Riemann zeta function.

The work of Bost and Connes sparked a resurgence in the study of KMS states for C*-dynamical systems constructed from algebraic and combinatorial structures, such as Bost--Connes type systems, Hecke systems associated to ring inclusions, and Toeplitz-type  systems associated to $ax+b$-semigroups, graphs, higher-rank graphs, C*-correspondences, among others. Most of these systems exhibit symmetries or various regularities that lead to a single natural choice for partition function. However, in some situations currently being considered this is no longer the case. The possibility of multiple partition functions was  recognized early on, because for systems containing isometries that are eigenvectors for the dynamics with nontrivial eigenvalues, the time evolutions are not approximately inner by \cite[Theorem 1]{OP} and there is no reason to expect a unique intrinsic Hamiltonian. Indeed, the evidence is that one should not even expect the energy levels to be uniquely determined. 
Because of this lack of uniqueness, it becomes necessary to have a more precise definition of what is really meant by the partition function of a general C*-dynamical system. Since in principle such systems can range from having no equilibrium states at all, to having multiple equilibrium states with very different Hilbert space representations, it is not surprising that the partition functions depend on the  states themselves. 

The equilibrium temperature space of a C*-dynamical system and, more specifically, the family of simplices of KMS$_\beta$ states, indexed by the inverse temperature $\beta \in [-\infty, \infty]$ can be quite arbitrary \cite{BEH,BEK}, and it is not clear what can be said in general about energy levels and partition functions. In this paper we advance the idea that if one focuses on the Murray--von Neumann type of the factor representations associated to extremal KMS states through the GNS construction, then certain invariants emerge that are both interesting and computable. Specifically, as a first step along the lines set out in \cite[Ch. 4]{CM} we concentrate on the extremal KMS$_\beta$ states whose GNS representations are factors of type I. The pure phases associated to such type I  states, although nonphysical, still behave in a fashion that agrees with our intuition at the macroscopic scale: when a C*-algebraic system is in one of these equilibrium states and the temperature is lowered, then equilibrium becomes more stable and the pure phase has continuity properties that link it to a unique ground state. 
It is for these type I factor KMS states that we define partition functions. The resulting families of  functions, fibered over the space of extremal type I states of the KMS simplex, constitute our basic C*-dynamical invariant for a C*-dynamical system. They are explicitly computable in terms of appropriate irreducible representations, and, of course, they are also invariant under isomorphisms of C*-dynamical systems, that is to say, under $\RR$-equivariant isomorphisms of the underlying C*-algebras.
 
 Partition functions as C*-dynamical invariants for systems arising from number theory have been considered before, e.g. \cite{CorMar}, but they have been largely used only for systems that have unique Hamiltonians due to the action of large symmetry groups, e.g.  Bost--Connes systems. Our present interest in partition functions originates in the study of the time evolution determined by the absolute norm on the C*-algebra of the $ax+b$-semigroup where $a$ and $b$ are algebraic integers in a given number field and $a$ is restricted to be in a given congruence monoid. These systems have been  introduced in \cite{Bru1} and  include the Toeplitz type systems from \cite{CDL} as a particular case. They exhibit a low-temperature phase transition \cite{Bru2} parametrized by the tracial states on certain group C*-algebras, which, intriguingly, also play a central role in the K-theory computation in \cite{BruLi}; reinforcing the connection made in the discussion following \cite[Theorem~6.6.1]{CELY} between the K-theory computations from \cite{CEL1,CEL2} for the case of the full $ax+b$-semigroup and the phase transition described in \cite{CDL}. Here we focus on the extremal KMS$_\beta$ states for $\beta >2$ and classify them according to type I or II, thus completing the type classification of KMS states started in \cite[Theorem~3.2(ii)]{Bru2}, see also  \cite[Remark~3.3(g)]{Bru2}. There are infinitely many extremal KMS$_\beta$ states of type I, and they give rise to an infinite family of partition functions  leading to various explicitly computable C*-dynamical invariants for congruence monoids and for the underlying number fields. 

We discuss next the organization of the paper and summarize its content, explaining roughly what is done in each section and what the main results are. 
In \secref{sec:partitionfunctions} we lay out what we mean by the partition function associated to a 
type I factor KMS state of a general C*-dynamical system, and we show how it can be computed in terms of an \emph{admissible triple} $(\H, \pi,\rho)$ consisting of a Hilbert space $\H$ on which there is an irreducible representation $\pi$ with a density operator $\rho$
for the state.  In  \proref{prop:gibbsbydensity} we show that for each type I factor KMS  state there is only one admissible triple up to unitary equivalence, and in \proref{KMSforsupercriticalbeta} we show that the KMS states are naturally grouped into quasi-equivalence classes that share a normalized Hamiltonian.  Taking limits as $\beta \to \infty$ in \corref{cor:kmsinfinityandground} we exhibit the KMS$_\infty$ states and we show that every admissible triple arises from the GNS representation of a pure ground state in one of these quasi-equivalence classes. Several families of examples are discussed in \secref{sec:examples} to illustrate the different possibilities for  partition functions.

In \secref{sec:congruencemonoids} we recall   the definition of congruence monoids $R_{\m,\Gamma}$ and the C*-algebras associated with their actions on rings of algebraic integers in \cite{Bru1}. The KMS$_\beta$ states of these systems have been characterized in \cite{Bru2}; for each $\beta>2$ they are indexed by the tracial states of a finite direct sum of group C*-algebras of the form $C^*(\a_\k\rtimes R_{\m,\Gamma}^*)$, where the $\a_\k$ are ideals representing ideal classes in a quotient of the ray class group modulo $\m$  and 
$R_{\m,\Gamma}^*$ is the group of units in the congruence monoid  $R_{\m,\Gamma}$, acting by multiplication.
We develop formulas for the KMS$_\beta$ states along the lines of those obtained in \cite{LR,CDL}. 
 These are needed later in order to identify convenient admissible triples  for the computation of the associated partition functions. The extremal KMS$_\beta$ states are indexed by extremal tracial states of the C*-algebras $C^*(\a_\k\rtimes R_{\m,\Gamma}^*)$.  In \secref{sec:traces} we discuss tracial states of direct sums and  in \thmref{thm:tracesandmeasures} we characterize the extremal KMS states of our systems in terms of ergodic invariant measures on tori, analogously to the characterization given in \cite{LaWa} for the systems from \cite{CDL}.

In \secref{sec:type} we construct concrete representations for the C*-algebras from actions of congruence monoids on rings of integers, \proref{prop:extendrep}, and  then parametrize the low-temperature KMS states expressly in terms of traces. The main technical results are
\thmref{thm:vNalg} and \proref{prop:tracestostates} which allow us to classify extremal KMS states according to type, \corref{cor:typeKMS}. Through the general theory developed in \secref{sec:partitionfunctions}, this also allows us to compute the partition functions of type I factor KMS states, \thmref{thm:partitionfncI}. These are the main  invariants for our C*-dynamical systems.  

We begin \secref{sec:topstructureKMS} with a discussion of the topological structure of the set of fixed points of an action. We then show that although the actions of the unit group $R_{\m,\Gamma}^*$ by  linear toral automorphisms of the (duals of) ideals in different ideal classes are not conjugate, there is nevertheless solidarity in terms of number of fixed points. By focusing on poles of the partition functions at critical $\beta$ we are able to select the subcollection consisting of the type I factor KMS states that have minimal residue. Our main result is \thmref{thm:main}, where we use this subset of minimal states and their partition functions to retrieve number-theoretic information about the initial congruence monoid and its underlying number field,  and we list concrete invariants that can be derived directly from the C*-dynamical invariants in \corref{cor:invariants}. 

As an application of our main results, we show in \secref{sec:recon} how several features of number fields can be characterized in terms of C*-dynamical invariants. The main application is \thmref{thm:recon}, where we obtain the Dedekind zeta function of the field, the generalized ideal class number, and the Kronecker set of a class field obtained naturally from the initial data defining the congruence monoid.
This allows us to characterize the field itself using C*-dynamical invariants.  That these quantities are invariant has been shown via K-theory, without any assumptions about time evolutions, in \cite{BruLi}. Here we reiterate the point that KMS states very often lead to results quite similar to those obtained via K-theory, through methods that are quite different, and that some of the number-theoretic invariants are computable directly from the representations associated to KMS states. These include the generalized class number, the zeta function of the modulus of the monoid, the partial zeta function of the trivial ideal class, and, up to finite sets of primes, the  Kronecker set of the class field associated with the initial data defining the congruence monoid.

\subsection*{Acknowledgment.} Part of this research was carried out during a visit of Bruce to Kyoto, and he would like to acknowledge this and to thank the faculty of Arts and Sciences at the Kyoto Institute of Technology for their hospitality and support. In addition, part of this research was performed during two visits of Takeishi to Victoria, and he would like to acknowledge the support and the hospitality of the Department of Mathematics and Statistics at the University of Victoria. 
The authors would also like to thank the referees for several very useful comments and suggestions.
\section{Persistence of Type I equilibrium under cooling and partition functions} 
\label{sec:partitionfunctions}
A \emph{time evolution} or {\em dynamics} on a C*-algebra $A$ is a group homomorphism $\sigma:\RR\to\Aut(A)$ such that for each fixed $a\in A$ the map  $t\mapsto\sigma_t(a)$ is norm continuous. The pair $(A,\sigma)$  is known as a \emph{C*-dynamical system}. The standard notion of equilibrium for these systems is given by KMS states. We refer the reader to \cite{BR} for general background on KMS states. Every C*-dynamical system determines a function that assigns to each $\beta\in \RR$  the closed convex set $\textup{KMS}_\beta(A,\sigma)$
of KMS$_\beta$ states. We will assume throughout that $A$ is unital, in which case this closed convex set is actually a Choquet simplex. In general, $\beta \mapsto \textup{KMS}_\beta(A,\sigma)$ can be a fairly arbitrary function, see e.g. \cite{BEK}, but if we restrict our attention to states of type I, i.e., those having GNS representations of type I in the Murray--von Neumann classification, then the extremal ones have a property that is quite intuitive in classical macroscopic terms: equilibrium `improves' when the temperature is lowered.  In this section we formulate this phenomenon explicitly and give a precise definition of the partition function that is naturally associated to each quasi-equivalence class of type I factor KMS states for general C*-dynamical systems.  We begin with the following slight reformulation of \cite[Lemma~4.3]{CCM}. 

\begin{proposition}\label{prop:gibbsbydensity} 
Let $\sigma$ be a time evolution on the unital, separable C*-algebra $A$. Fix $\beta_0 \in (0,\infty)$ and suppose that $\varphi$ is an extremal $\sigma$-KMS$_{\beta_0}$ state of type $I$.  Then 
\begin{enumerate}
\item[(i)] 
there exists an irreducible representation $\pi$ of $A$ on a separable Hilbert space $\H$ and a positive trace-class operator $\rho$ of norm $1$ on $\H$ such that 
\begin{equation}\label{eqn:gibbsdensity}
\varphi(a)=\frac{\Tr(\pi(a)\rho)}{\Tr(\rho)} \qquad a\in A;
\end{equation}
\item[(ii)]  
the operator $\rho$ is strictly positive and $U_t:= \rho^{-it/\beta_0}$ is a strongly continuous one-parameter unitary group  on $\H$ such that
\begin{equation}\label{eqn:timespatial}
  \pi(\sigma_t(a)) = U_t \pi(a) U_t^* \qquad a\in A, \  t\in\RR;
\end{equation}
\item[(iii)] any two triples $(\pi_j, \H_j, \rho_j)$ for $j =1,2$  arising as in part (i) above are unitarily equivalent, and this unitary
equivalence intertwines the unitary groups $\rho_1^{-it/\beta_0}$ and $\rho_2^{-it/\beta_0}$.
\end{enumerate}
\end{proposition}
\begin{proof}
Let $(\pi_\varphi, \H_\varphi, \xi_\varphi)$ be the  GNS representation and cyclic vector associated to $\varphi$.
By assumption, $\M_\varphi:=\pi_\varphi(A)''\subseteq \B(\H_\varphi)$ is a type I factor, so there exists an irreducible representation
 $\pi$ of $A$ on a separable Hilbert space $\H$ such that the map $\theta: \pi_\varphi(a) \mapsto \pi(a)$ extends to a von Neumann algebra isomorphism $\theta: M_\varphi \to \B(\H)$. 
 
 The  vector state 
 $\tilde{\varphi} (x) = \langle x\xi_\varphi, \xi_\varphi\rangle$ is a  normal state of $\M_\varphi = \pi_\varphi(A)''$. Hence $\tilde{\varphi}\circ \theta\inv$ is a  normal state of $\B(\H)$ so there exists a  positive trace class operator $\rho$ on $\H$, such that
 \begin{equation}\label{densityformula}
\tilde\varphi  (\theta\inv ( x) ) = \frac{\Tr( x\rho)}{\Tr(\rho)}.
\end{equation}
Writing $x = \pi(a)$ above and noticing that $\tilde\varphi(\theta\inv(\pi(a))) = \langle \pi_\varphi(a)\xi_\varphi, \xi_\varphi\rangle =  \varphi (a)$ we see that \eqref{eqn:gibbsdensity} holds for every $a\in A$. We point out that \eqref{eqn:gibbsdensity} only determines 
the operator $\rho$ up to a scalar multiple, and we choose to normalize it by requiring $\|\rho\| =1$. 
This completes the proof of part (i).

 By \cite[Corollary~5.3.9]{BR}, the cyclic vector $\xi_\varphi$ is separating for $\M_\varphi$, so 
 $\tilde{\varphi} (x) = \langle x\xi_\varphi, \xi_\varphi\rangle$ is a faithful state of $\M_\varphi$ and 
 $\tilde{\varphi}\circ \theta\inv$ is a  faithful normal state of $\B(\H)$. 
Thus $\rho$ is a strictly positive compact operator of norm 1 on $\H$, i.e.
there exists an orthonormal basis $\{\xi_n\}_{n\in \NN}$ of $\H$ consisting of 
eigenvectors with eigenvalues $\omega_n \in (0,1]$.
It follows that for each complex number $z$ with $\Im z \geq 0$ the map $U_z: \xi_n \mapsto {\omega_n^{-iz/\beta_0}} \xi_n$ extends by linearity and continuity to a contraction 
 $U_z$ on $\H$.  
 
A standard approximation argument then shows that $\{U_z\}_{\Re z \geq 0}$  is a strongly continuous contraction semigroup  on $\H$ whose restriction to the real axis is a unitary group, and whose restriction to the open upper half-plane is an analytic function. Moreover,  $\rho = U_{i\beta_0}$ so the state of $\B(\H)$ defined by \eqref{densityformula}  is of Gibbs type with respect to the unitary group $U$ and hence satisfies the KMS$_{\beta_0}$ condition for the associated time evolution $\Ad_{U_t}$ on $\B(\H)$, which is sigma-weakly continuous. This can be verified directly: given $X$ and $Y$ in $\B(\H)$, the  analytic function interpolating between the boundary values $F_{X,Y}(t+ i0) = \Tr(X U_t Y U_
t^* \rho)$ and 
$F_{X,Y}(t+ i\beta_0) = \Tr( U_t Y U_t^* X\rho)$ for $t\in \RR$ is given by $F_{X,Y}(z) = \Tr(X U_z Y U_{\overline z}^* \rho)$ . 

On the other hand, at the level of the GNS representation of $\varphi$, a standard argument gives a 
  unitary group $U^\varphi$ on $\H_\varphi$ defined by 
  $U^\varphi_t : \pi_\varphi(b) \xi_\varphi  \mapsto \pi_\varphi(\sigma_t( b)) \xi_\varphi$ for $b\in A$ and $t\in \RR$. The associated  time 
evolution $\tilde\sigma := \Ad_{U^\varphi_t}$ on $\M_\varphi$ is sigma-weakly continuous and satisfies $\tilde\sigma_t(\pi_\varphi(a)) =  \pi_\varphi(\sigma_t(a))$ 
and the state $\tilde\varphi$ is a $\tilde\sigma$-KMS$_{\beta_0}$ state of $\M_\varphi$, see e.g. \cite[Corollary 5.3.4]{BR}. 
Then $\tilde\varphi\circ \theta\inv$ is also a $(\theta \circ \tilde\sigma \circ\theta\inv)$-KMS$_{\beta_0}$ state on $\B(\H)$. 
By \cite[Theorem 13.2]{T}, (see also \cite[Theorem 5.3.10]{BR}), the time evolution is determined by the KMS$_{\beta_0}$ state, so  $\theta \circ \tilde\sigma\circ \theta\inv = \Ad_{U_t}$ on $\B(\H)$, from which  \eqref{eqn:timespatial} follows on restricting to $\pi(A)$ because 
$\theta \circ \tilde\sigma_t\circ \theta\inv (\pi(a)) = \theta \circ \tilde \sigma_t (\pi_\varphi(a)) = \theta (\pi_\varphi(\sigma_t(a))) = \pi(\sigma_t(a))$. This completes the proof of part (ii).

For part (iii), simply observe that  $\pi_1$ and $\pi_2$ are two irreducible representations that are quasi-equivalent to $\pi_\varphi$,
so they are mutually unitarily equivalent  via a unitary $\U :\H_1 \to \H_2$. Then $\U\rho_1\U^*$ is a strictly positive 
trace-class operator of norm 1 on $\H_2$ such that 
\[
 \frac{\Tr( \pi_2(a) \, \U\rho_1\U^*)}{\Tr(\U\rho_1\U^*)} =  \frac{\Tr(\U^* \pi_2(a) \, \U\rho_1)}{\Tr(\rho_1)} =  \frac{\Tr( \pi_1(a) \rho_1)}{\Tr(\rho_1)} = \varphi(a) = \frac{\Tr( \pi_2(a) \rho_2)}{\Tr(\rho_2)}
\]
for every $a\in A$. 
Since $\pi_2(A)$ is dense in $\B(\H_2)$ and each normal state of $\B(\H_2)$ has a unique density operator,  
we have  $\frac{1}{\Tr(\rho_1)}\U\rho_1\U^* =  \frac{1}{\Tr(\rho_2)}\rho_2$, and since both $\U\rho_1\U^*$ and $\rho_2$ have norm $1$, they are equal.
\end{proof}

\begin{remark}
 The result and its proof are essentially from \cite[Lemma~4.3]{CCM}, see also \cite[Proposition~4.52]{CM}. Here we have taken into account that the regularity property \cite[Definition~4.1]{CCM} is automatic for KMS states  \cite[Corollary~5.3.9]{BR} and  we have stated explicitly the uniqueness property of the construction. We have also removed the assumption that for sufficiently large $\beta $ {\em all} extremal KMS$_\beta$ states are of type I. This assumption is stated shortly before \cite[Lemma~4.3]{CCM}, see also \cite[Remark 4.50]{CM}, but it is not necessary here, which is important in our applications.
\end{remark}
\begin{remark}
If the C*-algebra $A$ has an extremal $\sigma$-KMS$_{\beta_0}$ state of type I for some $\beta_0 <0$, then \proref{prop:gibbsbydensity} can be applied to the reverse dynamics $\sigma_{-t}$.
\end{remark}

Fixing the notation $(\pi, \H, \rho)$ to denote the admissible triple associated to a KMS$_{\beta_0}$ state $\varphi$
in \proref{prop:gibbsbydensity}, 
the infinitesimal generator of the unitary group $U_t = \rho^{-it/\beta_0}$ is the operator
\begin{equation} \label{HfromOmega}
H:= {\textstyle\frac{-1}{\beta_0}}\log \rho
\end{equation}
and we may write the density operator and the unitary group in terms of $H$ as $\rho = e^{-\beta_0 H} $ and  
$\rho^{-it/\beta_0} = e^{itH} $, respectively. Because of this we view $H$ as a `represented Hamiltonian'
associated to the state $\varphi$ that can be regarded as the `stable half' of the Liouville operator $H_\varphi$ 
studied in \cite[Chap. V]{tB}, see Proposition~\ref{pro:gnsfactorization} below.

The normalization $\|\rho\| =1$ is somewhat unusual, but it serves the purpose of
forcing all the represented Hamiltonians that arise as in Proposition~\ref{prop:gibbsbydensity}
from extremal KMS$_\beta$ states of type I  
to have $0$ as the minimum of their spectra. 
In general, if $e^{-\beta H}$ is of trace class and $\beta' \geq \beta$, then $e^{-\beta' H}$ is also of trace class, so
we define $\beta_c: = \inf \{ \beta\in [-\infty, \infty): \Tr(e^{-\beta H}) <\infty\} $. Then $0\leq \beta_c \leq \beta_0$, unless the Hilbert space is finite-dimensional, in which case $\beta_c = -\infty$, so that $\beta_c \in [0,\infty) \cup \{-\infty\}$.
We shall see next that there is a family  consisting of exactly one KMS$_\beta$ state for each positive $\beta \in (\beta_c, \infty) $
that share the same Hamiltonian with $\varphi$ (up to unitary equivalence). 

\begin{definition}
\label{def:partitionfunction}
Let $\varphi$ be an extremal $\sigma$-KMS$_{\beta_0}$ state of type I, and let $H$ be the represented Hamiltonian associated to $\varphi$ in Equation~\eqref{HfromOmega}. The {\em partition function of $\varphi$} is the function $Z_\varphi(s):= \Tr(e^{-s H})$ defined for $\Re s > \beta_c$.
\end{definition}

\begin{proposition}\label{KMSforsupercriticalbeta}
Fix $\beta_0 \in (0,\infty)$ and suppose \tcb{$\varphi$} is an extremal KMS$_{\beta_0}$ state of type I of the system $(A,\sigma)$,   let  \tcb{$(\pi,\H, \rho)$} be the triple associated to \tcb{$\varphi$} in Proposition~\ref{prop:gibbsbydensity}, and let \tcb{$H:= \frac{-1}{\beta_0}\log \rho$} be the \tcb{represented}  Hamiltonian. 
For each positive $\beta$  such that \tcb{$e^{-\beta H} (= \rho^{\beta/\beta_0})$} is of trace-class 
define a state on $A$ by 
\begin{equation}\label{cooledstate}
\varphi_{\beta_0, \beta} (a)  := \tcb{\frac{\Tr( \pi(a) e^{-\beta H})}{\Tr (e^{-\beta H})} }\qquad a \in A.
\end{equation}
Then $\varphi_{\beta_0, \beta}$ is an extremal $\sigma$-KMS$_{\beta}$ state of type I and 
  \tcb{$(\pi,\H, \rho^{\beta/\beta_0})$} is an admissible triple for $\varphi_{\beta_0, \beta}$ in the sense of Proposition~\ref{prop:gibbsbydensity}
so all the states $\varphi_{\beta_0, \beta}$ have \tcb{$H$} as represented Hamiltonian
and $Z_{\varphi_{\beta_0, \beta}}(s)=\tcb{\Tr (e^{-s H})}$ for every $\beta >\beta_c$.
Moreover, the quasi-equivalence class of the state \tcb{$\varphi$} among all KMS$_\beta$ states for $\beta>0$ is precisely $\{ \varphi_{\beta_0, \beta}: \beta >\beta_c\}$.
\end{proposition}

\begin{proof} Formula \eqref{cooledstate} clearly defines a state that is quasi-equivalent to \tcb{$\varphi$} because  it is given by a density operator in the same irreducible representation, and hence is a type I factor state. A computation analogous to the one in the proof of Proposition~\ref{prop:gibbsbydensity} shows that this state satisfies the KMS$_\beta$ condition for the time evolution $\sigma$ implemented by \tcb{$e^{itH}$}. 
Thus, \tcb{$(\pi, \H, e^{-\beta H})$} is an admissible triple for $\varphi_{\beta_0, \beta}$ and its associated Hamiltonian
is \tcb{$ {\textstyle\frac{-1}{\beta}}\log e^{-\beta H} = H$}.
To prove that the quasi-equivalence class of \tcb{$\varphi$} is the given set, assume $\psi$ is a KMS$_\beta$ state of type I for some $\beta>0$ and that $\psi$ is quasi-equivalent to \tcb{$\varphi$}. Then $\psi$ is a factor state, so it is extremal and it follows from \cite[Theorem 5.3.30(4)]{BR} that $\psi_{\beta,\beta'} = \varphi_{\beta_0, \beta'}$ for all $\beta' > \max\{\beta,\beta_c\}$. Hence, $\psi$ and \tcb{$\varphi$} share the same Hamiltonian, which implies that $\psi  =\varphi_{\beta_0, \beta}$ and that \tcb{$e^{-\beta H}$} is trace class.
\end{proof}
\begin{remark}
The partition function of an extremal $\sigma$-KMS$_{\beta_0}$ state $\varphi$ of type I from Definition~\ref{def:partitionfunction} coincides with the one mentioned in \cite[Remark~2.2(ii)]{LLN}. We emphasize here that   the operator $H$ and hence the partition function $Z_\varphi$ depend only on the quasi-equivalence class of $\varphi$, which consists of the states that have densities $\frac{e^{-\beta H}}{\Tr(e^{-\beta H})}$  on the Hilbert space of $\varphi$,
cf. Proposition~\ref{KMSforsupercriticalbeta}.
Thus, to obtain $Z_\varphi$, it suffices to find a triple $(\pi,\H,H)$ that is admissible for $\varphi$ as in Proposition~\ref{prop:gibbsbydensity}.
This is useful when one wishes to make explicit computations.
\end{remark}

By definition, following \cite{CM},  we reserve the terminology {\em KMS$_\infty$ states} to refer to the weak-* limits of KMS$_\beta$ states as $\beta \to \infty$. Next we extend Proposition~\ref{KMSforsupercriticalbeta} to cover the case $\beta = \infty$ by interpreting it as a such a limit.

\begin{corollary}\label{cor:kmsinfinityandground}
In the conditions of Proposition~\ref{KMSforsupercriticalbeta}, let $Q_0$ be the  spectral projection of the largest eigenvalue $\lambda_0 =1$ of $\rho$, and let $\xi$ be a unit vector in the range of $Q_0$. Then  
\begin{enumerate}
\item the state $\varphi_{\beta_0, \infty}$ on $A$ defined by 
\begin{equation}\label{kmsinfstate}
\varphi_{\beta_0, \infty} (a)  := \frac{\Tr( \tcb{\pi}(a) Q_0)}{\Tr (Q_0)} \qquad a \in A
\end{equation}
is a type I factor KMS$_\infty$ state  that is quasi-equivalent to $\varphi_{\beta_0,\beta}$ and 
$\lim_{\beta \to \infty} \varphi_{\beta_0,\beta} = \varphi_{\beta_0, \infty}$; 
\item the vector state $\omega_\xi$ defined by 
\[\omega_{\xi}(a) := \langle \tcb{\pi}(a) \xi, \xi\rangle  \qquad a \in A
\]
is a pure ground state, and $\tcb{\pi}$ together with the cyclic vector $\xi$ is unitarily equivalent to the GNS representation of $\omega_\xi$.
\end{enumerate}
\end{corollary}
\begin{proof}
Since $\varphi_{\beta_0, \infty}$ is defined by a density in an irreducible representation, it is a type I factor state  and is quasi-equivalent to $ \varphi_{\beta_0,\beta}$. 
The distance 
$
\|\varphi_{\beta_0, \beta} - \varphi_{\beta_0, \infty}\| $ between the states  
is given by the trace-norm of the difference of densities. Since
\[
 \left\|  \frac{ e^{-\beta H}}{\Tr (e^{-\beta H})} - \frac{ Q_0}{\Tr (Q_0)} \right\|_1 
\leq    \left|\frac{ 1}{\Tr (e^{-\beta H})} - \frac{1}{\Tr (Q_0)} \right| \|e^{-\beta H}\|_1  + \frac{1}{\Tr (Q_0)} \| e^{-\beta H} - Q_0   \|_1
\]
and since  $\| e^{-\beta H} - Q_0   \|_1 = \Tr (e^{-\beta H}) - \Tr (Q_0)$, in order to show that 
$
\|\varphi_{\beta_0, \beta} - \varphi_{\beta_0, \infty}\| \to 0 $, and hence that $\phi_{\beta,\infty}$ is a KMS$_\infty$ state,
 it suffices to show that 
$\lim_{\beta \to \infty} \Tr (e^{-\beta H}) = \Tr (Q_0)$.
To do this, let  $\{\xi_k\}_{k\in \NN}$ be an orthonormal basis of $\H$ diagonalizing $H=\textup{diag}(d_0 =0, d_1, d_2,...)$ with $d_k := \frac{-\log \lambda_k}{\beta_0}$ strictly increasing and let $a_k$ be the multiplicity of $d_k$. Suppose $\beta \geq \beta_0$. Then
\[
0 \leq \Tr (e^{-\beta H}) - \Tr (Q_0) = Z_\varphi(\beta) - a_0 =\sum_{k=1}^\infty a_ke^{-\beta d_k} \leq 
 e^{-\beta d_1}\sum_{k=1}^\infty a_ke^{-\beta_0 (d_k-d_1)} \to 0 
\]
 because  $e^{-\beta d_1} \to 0$  as $\beta \to \infty$ and  $\sum_{k=1}^\infty a_ke^{-\beta_0 (d_k-d_1)}$ converges. This completes the proof of part (1).
 
To prove part (2) notice that  $\xi$ is cyclic  because  $\tcb{\pi}$ is irreducible so, by the essential uniqueness of the GNS representation of  $\omega_\xi$, it  is unitarily equivalent to $\tcb{\pi}$, hence $\omega_\xi$ is pure. It remains to show that $\omega_\xi$ is a ground state.
Observe first that  
$\varphi_{\beta_0, \infty} = \frac{1}{\Tr(Q_0)}\sum_i \omega_{\xi_i}$, where the sum is over an orthonormal basis of the range of $Q_0$. Then $\omega_\xi$ is dominated by a positive multiple of $\varphi_{\beta_0, \infty}$, which is a KMS$_\infty$ state and hence a ground state itself. Since the set of ground states in a unital algebra is a face of the set of states \cite[Theorem 5.3.37]{BR}, it follows that $\omega_{\xi}$ is a ground state.
\end{proof}

\begin{proposition}\label{pro:gnsfactorization}
Suppose that $(\pi,\H,\rho)$ is a triple arising from an extremal $\sigma$-KMS$_\beta$ state $\varphi$ of type I, as in Proposition~\ref{prop:gibbsbydensity} with $\rho= e^{-\beta H}$. For $t\in\RR$, let $W_t:=e^{itH}\otimes e^{-itH}$, so that $\{W_t\}_{t\in \RR}$ is a continuous one-parameter unitary group on $\H\otimes\H$, let  $\{\xi_0, \xi_1, \ldots , \xi_n, \ldots\}$  
be an orthonormal basis for $\H$ consisting of eigenvectors of $H$
 with eigenvalues   $\{0 = h_0 \leq h_1 \leq \cdots \leq h_n \leq \cdots\}$ and 
let 
\[
\eta:=\frac{1}{\Tr(e^{-\beta H})^{1/2}}\sum_{n=0}^\infty e^{\frac{-\beta h_n}{2}}\xi_n\otimes\xi_n.
\]
Then there exists a unitary $\mathcal{V}:\H_\varphi\overset{\cong}{\to}\H\otimes\H$ such that 
\begin{enumerate}
\item[(i)] $\mathcal{V}\xi_\varphi=\eta$;
\item[(ii)] $\mathcal{V}\pi_\varphi(a)\mathcal{V}^*=\pi(a)\otimes 1$ for all $a\in A$;
\item[(iii)] $\mathcal{V}U_t^\varphi\mathcal{V}^*=W_t$ for all $t\in\RR$.
\end{enumerate}
\end{proposition}

\begin{proof}
Let $\Phi:\B(\H)\to\B(\H\otimes\H)$ by $\Phi(T)=T\otimes 1$. Then
\begin{align*}
\langle\Phi\circ\pi(a)\eta,\eta\rangle&=\frac{\Tr(\pi(a)e^{-\beta H})}{\Tr(e^{-\beta H})}\\
&= \varphi(a) \quad\text{by \eqref{eqn:gibbsdensity}.}
\end{align*}
By assumption $\pi(A)''=\B(\H)$, hence $(\Phi\circ \pi(A))'' = \B(\H)\otimes 1$. It is easy to see that for any pair $(m,n)$ there exists a rank-one operator $T \in B(\H)$  such that $(T \otimes 1)\eta = \xi_m \otimes \xi_n$, so  the unit vector $\eta\in\H\otimes\H$ is cyclic for the representation $\Phi\circ\pi$. By uniqueness of the GNS representation, there is a unitary $\mathcal{V}:\H_\varphi\overset{\cong}{\to}\H\otimes\H$ such that (i) and (ii) hold. The following calculation shows that (iii) holds:
\begin{align*}
W_t\mathcal{V}(\pi_\varphi(a)\xi_\varphi)&=W_t(\pi(a)\otimes 1)\eta\\
&=W_t\left(\frac{1}{\Tr(e^{-\beta  H})^{1/2}}\sum_{n=0}^\infty e^{\frac{-\beta h_n}{2}}(\pi(a)\xi_n)\otimes \xi_n\right)\\
&=\frac{1}{\Tr(e^{-\beta  H})^{1/2}}\sum_{n=0}^\infty e^{\frac{-\beta h_n}{2}}(e^{itH}\pi(a)\xi_n)\otimes (e^{-iH }\xi_n)\\
&=\frac{1}{\Tr(e^{-\beta  H})^{1/2}}\sum_{n=0}^\infty e^{\frac{-\beta h_n}{2}}(\pi(\sigma_t(a))e^{ith_n}\xi_n)\otimes (e^{-ith_n }\xi_n) \\
&=\frac{1}{\Tr(e^{-\beta  H})^{1/2}}\sum_{n=0}^\infty e^{\frac{-\beta h_n}{2}}(\pi(\sigma_t(a))\xi_n)\otimes \xi_n \\
&=\big(\pi( \sigma_t( a )) \otimes 1\big)\eta\\
&=\big(\mathcal{V}\pi_\varphi(\sigma_t(a))\mathcal{V}^*\big)\eta\\
&=\mathcal{V}U^\phi_t(\pi_\varphi(a))\xi_\varphi.
\end{align*}
\end{proof}

\begin{remark} Recall from \cite{tB,tBW} that if $\varphi$ is an invariant state of a C*-dynamical system $(A,\sigma)$, then  the infinitesimal generator of the usual unitary group defined by
$U^\varphi_t : \pi_\varphi(b) \xi_\varphi  \mapsto \pi_\varphi(\sigma_t( b)) \xi_\varphi$ for $b\in A$ and $t\in \RR$
on the GNS Hilbert space of $\varphi$  is called the {\em Liouville operator} associated to $\varphi$.
By \proref{pro:gnsfactorization}  the spectrum of the Liouville operator associated to $\varphi$ 
is given by $\{h_j -h_k: j,k  = 0, 1, 2, \ldots\}$, and is symmetric with respect to the origin, cf. \cite{tBW}.
The spectra of the represented admissible Hamiltonians for $\varphi$ are the same for quasi-equivalent type I factor KMS states, and do not  depend on the inverse temperature above critical; thus, the same is true of the spectra of the associated Liouville operators, cf. \cite[Theorem A]{tBW}. However, as we shall see in \secref{sec:type}, for the examples arising from congruence monoids, the represented Hamiltonians  of two inequivalent type I factor KMS states  
do not have to have the same spectral sets. We suspect the same may be true for the spectra of  the associated Liouville operators, even though the corresponding representations are faithful by \thmref{thm:partitionfncI} below. This would not contradict \cite[Theorem A]{tBW} because our systems are not approximately inner by \cite[Theorem 1]{OP}.
\end{remark}

\begin{remark}
A consequence of \corref{cor:kmsinfinityandground} is that the admissible triples considered in \proref{prop:gibbsbydensity} always arise as the GNS representations of those pure ground states that are quasi-equivalent to a KMS$_\beta$ state for finite $\beta$. This shows that the admissible triples  are \tcb{intrinsic} to the C*-dynamical system in the sense that they can be constructed from 
any ground state that lies below a KMS$_\infty$ state.
In general, not every ground state gives rise to an admissible triple.
To see why, consider, for instance,   $\O_\infty$ with the usual periodic dynamics, in which there is a ground state but no KMS$_\beta$ states (see \cite{OP}).
\end{remark}

\section{Examples of partition functions}
\label{sec:examples}
{\bf Semigroup C*-algebras of quasi-lattice ordered monoids.}
Let $(G,P)$ be a quasi-lattice ordered group, and \tcb{let $C_\lambda^*(P)$ denote the left regular C*-algebra of $P$, which by definition is the C*-algebra generated by the image of the  left regular representation $\lambda$ of $P$ on $\ell^2(P)$. If $N:P\to [1,\infty)$ is a multiplicative map such that $N(p)=1$ only if $p=e$, then there is an associated dynamics  $\alpha^N$  on $C_\lambda^*(P)$ determined by  $\alpha^N( \lambda_p) = N(p)^{it} \lambda_p$. }Assume that the \emph{growth series} $Z_N(\beta):=\sum_{p\in P}N(p)^{-\beta}$ has abscissa of convergence $\beta_c$ for some $0<\beta_c<\infty$ (see \cite{BLRS} for details). For $\beta_0\in (\beta_c,\infty)$, let $\psi_{\beta_0}$ be the unique $\alpha^N$-KMS$_{\beta_0}$ state on $C_\lambda^*(P)$ from \cite[Definition~3.4]{BLRS}. Explicitly, $\psi_{\beta_0}(\cdot)=\frac{1}{Z_N(\beta)}\Tr(\ \cdot \ e^{-\beta H_N})$ where $H_N$ is the diagonal operator on $\ell^2(P)$ such that $H_N\delta_p=(\log N(p))\delta_p$ for all $p\in P$.

Using \cite[Theorem~3.5(1)]{BLRS}, it is easy to see that the left regular \tcb{representation  gives} rise to a triple $(\lambda,\ell^2(P), e^{-\beta H_N})$ satisfying Proposition~\ref{prop:gibbsbydensity}(i) and that the partition function is $Z_{\psi_{\beta_0}}(\beta)=Z_N(\beta)$, and is defined for all $\beta>\beta_c$. See \cite[\S~4]{BLRS} for more on the growth series $Z_N(\beta)$, including its relation to cliques in the special case that $P$ is a finitely generated right-angled Artin monoid.

{\bf The Bost--Connes system for $\QQ$.}
For each rational prime $p$, let $\ZZ_p$ denote the ring of $p$-adic integers, and let $\hat{\ZZ}:=\prod_p\ZZ_p$ be the ring of integral adeles over $\QQ$. One model for the Bost--Connes C*-algebra from \cite{BC} is the semigroup crossed product $C(\hat{\ZZ})\rtimes \NN^\times$ where $\NN^\times:=\NN\setminus\{0\}$ acts on $\hat{\ZZ}$ by multiplication (see \cite[\S~5.4]{La98} and \cite[\S~3]{LR99}).
For each $n\in\NN^\times$, let $\mu_n$ denote the corresponding isometry in $C(\hat{\ZZ})\rtimes \NN^\times$. The C*-algebra $C(\hat{\ZZ})\rtimes \NN^\times$ carries a canonical time evolution $\sigma^\QQ$ such that $\sigma^\QQ_t(\mu_n)=n^{it}\mu_n$ for all $n\in\NN^\times$ and $\sigma^\QQ_t(f)=f$ for all $f\in C(\hat{\ZZ})$. 
Assume $\beta>1$. Then the set of extremal points of the simplex of KMS$_\beta$ states on $C(\hat{\ZZ})\rtimes \NN^\times$ is homeomorphic to $\hat{\ZZ}^*:=\prod_p\ZZ_p^*$ where $\ZZ_p^*$ denotes the group of units in $\hat{\ZZ}_p$ (see \cite[Theorem~5]{BC}). Explicitly, given $\mathbf{u}\in\hat{\ZZ}^*$, the corresponding KMS$_\beta$ state $\omega_{\beta,\mathbf{u}}$ can be realized as follows (see \cite[Theorem~25]{BC}).
There is a representation $\pi_\mathbf{u}$ of $C(\hat{\ZZ})\rtimes \NN^\times$ on $\ell^2(\NN^\times)$ determined on the canonical orthonormal basis by $\pi_\mathbf{u}(f)\delta_k=f(k\mathbf{u})\delta_k$ for $f\in C(\hat{\ZZ})$ and $\pi_\mathbf{u}(\mu_n)\delta_k=\delta_{nk}$ for $n\in\NN^\times$. Let $H$ denote the diagonal operator $H\delta_n=(\log n)\delta_n$, so that $\zeta(\beta)=\Tr(e^{-\beta H})$ for $\beta>1$ where $\zeta(s)$ denotes the Riemann zeta function, and define
\[
\omega_{\beta,\mathbf{u}}(\cdot):=\frac{1}{\zeta(\beta)}\Tr(\pi_\mathbf{u}(\cdot)e^{-\beta H}).
\]
Then $(\pi_\mathbf{u},\ell^2(\NN^\times),e^{-\beta H})$ is a triple satisfying (i) from Proposition~\ref{prop:gibbsbydensity}, so the partition function of $\omega_{\beta,\mathbf{u}}$ is given by $Z_{\omega_{\beta,\mathbf{u}}}(\beta)=\zeta(\beta)$ for every $\mathbf{u}\in\hat{\ZZ}^*$.

We have restricted our attention here to $\QQ$ only for ease of exposition and accessibility. For the general Bost--Connes type system associated with a number field $K$ (as defined in \cite{HP,LLN}), one can similarly show that the partition function of any low-temperature KMS state is $\zeta_K(s)$, the Dedekind zeta function of $K$ (cf. \cite[Remark~2.2(ii)]{LLN}).

{\bf The left regular C*-algebra of $\NN\rtimes\NN^\times$.}
The left regular C*-algebra $C_\lambda^*(\NN\rtimes\NN^\times)$ carries a canonical time evolution $\sigma$ such that $\sigma(\lambda_{(m,k)})=k^{it}\lambda_{(m,k)}$ for all $(m,k)\in \NN\rtimes\NN^\times$ and $t\in\RR$. By \cite[Theorem~7.1(3)]{LR}, for each $\beta\in(2,\infty]$, the set of KMS$_\beta$ states of the system $(C_\lambda^*(\NN\rtimes\NN^\times),\sigma)$ is isomorphic to the simplex of probability measures on $\TT$, so that the extremal KMS$_\beta$ states are in bijection with the points in $\TT$, and the KMS$_\beta$ state $\phi_{\beta,z}$ corresponding to $z\in\TT$ is of type I. The state $\phi_{\beta,z}$ can be described as follows (see the proof of \cite[Proposition~9.3]{LR}).

Let $X=\{(r,x) : x\in \NN^\times, 0\leq r\leq x-1\}$ and let $\{\delta_{(r,x)} : (r,x)\in X\}$ be the canonical orthonormal basis for $\ell^2(X)$. For each $z\in\TT$, there is an irreducible representation $\pi_z$ of $C_\lambda^*(\NN\rtimes\NN^\times)$ on $\ell^2(X)$ such that 
\[
\pi_z(\lambda_{(1,1)})\delta_{(r,x)}=\begin{cases}
\delta_{(r+1,x)} & \text{ if } r<x-1,\\
z\delta_{(0,x)} & \text{ if } r=x-1,
\end{cases}
\] 
and
\[
\pi_z(\lambda_{(0,p)})\delta_{(r,x)}=\delta_{(pr,px)}
\]
for \tcb{every prime number $p$}. 
Let $H$ be the diagonal operator such that $H\delta_{(r,x)}=(\log x)\delta_{(r,x)}$, so that $\Tr(e^{-\beta H})=\zeta(\beta-1)$ where $\zeta(s)$ is the Riemann zeta function. Then 
\[
\phi_{\beta,z}(\cdot)=\frac{1}{\zeta(\beta-1)}\Tr(\pi_z(\cdot) e^{-\beta H}),
\] 
and $(\pi_z,\ell^2(X),e^{-\beta H})$ is a triple satisfying (i) from Proposition~\ref{prop:gibbsbydensity}, so the partition function of $\phi_{\beta,z}$ is given by $Z_{\phi_{\beta,z}}(\beta)=\zeta(\beta-1)$ for every $z\in \TT$.\\

{\bf The C*-algebra of the $ax+b$-semigroup over the ring of algebraic integers.}
Let $R$ be the ring of algebraic integers in a number field $K$. The left regular semigroup C*-algebra associated to the $ax+b$-semigroup $R\rtimes R^\times$ of the ring $R$ has been studied in \cite{CDL}.  \tcb{In the next section we will study in detail a generalization of this  considered in \cite{Bru1,Bru2}, so we will be brief on details here, 
but we would like to mention as motivation that the extremal KMS$_\beta$ states of $C_\lambda^*(R\rtimes R^\times)$ for high $\beta$ have been parametrized, first by extremal traces on the C*-algebra $\bigoplus_{\k \in\Cl} C^*(\a_\k \rtimes R^*)$ \cite[Theorem 7.3]{CDL} and subsequently by ergodic invariant probability measures for the action of  the unit group $R^*$ of $R$ on the tori $\widehat{\a}_\k$ together with characters of the associated isotropy subgroups \cite[Theorem 2.2]{LaWa}.  From the analysis in \cite{CDL} it is evident that the partial zeta functions associated to ideal classes are related to partition functions. Nevertheless, a concrete link has not been established formally until now, partly because, in contrast to previously studied systems, there is not a unique choice for partition function, but  several possibilities instead, that depend on the ideal class and the size of the orbit  associated to an extremal KMS state of type I.} In Theorem~\ref{thm:partitionfncI} below we will clarify the situation and describe the partition functions associated to all the extremal KMS states of type I for the systems introduced in \cite{Bru1,Bru2}, which contain, as special cases, those of \cite{CDL}.

\section{C*-algebras from actions of congruence monoids}\label{sec:congruencemonoids}
We begin with a brief review of some basic facts about congruence monoids and the semigroup C*-algebras associated with their actions on rings of algebraic integers \cite{Bru1}, and then we describe the phase transition of KMS states obtained in \cite{Bru2}.

Let $K$ be an algebraic number field with ring of integers $R$. A modulus $\m =  \m_\infty \m_0$ for $K$ consists of a $\{0,1\}$-valued function $\m_\infty$ on the set $V_{K,\RR}$ of real embeddings of $K$ together with a non-zero ideal $\m_0 = \prod_\p \p^{v_\p(\m_0)}$ in $R$, which can be viewed as the finitely supported function on prime ideals with integer values $\m_0(\p) = v_\p(\m_0)$.  As customary, the real embeddings in the support of $\m_\infty$ are said to divide $\m_\infty$, so we write $w\mid \m_\infty$ to mean $\m_\infty (w) = 1$.
  
Any $a \in R$ that is relatively prime to $\m_0$  can be {\em reduced modulo $\m$} by taking its sign in each real embedding that divides $\m_\infty$ together with its residue class modulo $\m_0$, that is, $[a]_\m := \left(\big(\sign w(a)\big)_{w\mid \m_\infty} , a+\m_0 \right)$ as an element of the multiplicative group of residues modulo $ \m$
\[
(R/\m)^* := \Big( \prod_{w \mid \m_\infty} \{ \pm1\}\Big) \times \big( R/\m_0\big)^*.
\] 
Given $\m$ and a subgroup $\Gamma$ of $\big( R/\m\big)^*$ we may implement a further restriction on 
residues modulo $\m$ by considering  the multiplicative submonoid of $R^\times$ consisting of elements that are relatively prime to $\m_0$ and reduce to an element of $\Gamma$ modulo $\m$; this is the {\em congruence monoid}
\[
R_{\m,\Gamma}:=\{ a \in R^\times: v_\p(a) =0 \text{ if } \p \mid \m_0 \text { and } [a]_\m \in \Gamma\}.
\]
The associated group of quotients  is 
\[
 R_{\m,\Gamma}\inv R_{\m,\Gamma} =K_{\m,\Gamma} :=\{ x \in K^\times: v_\p(x) =0 \text{ if }  \p \mid \m_0 \text { and } [x]_\m \in \Gamma\} .
 \]
See \cite[Ch.2.11]{GeHa} for the general background on congruence monoids and also \cite[Section 2.2]{Bru1} and \cite{BruLi} for the details of this construction in relation to C*-algebras.

Let $\I_\m$ be the group of fractional ideals in $K$ that are relatively prime to $\m_0$ and let $i(K_{\m,1})$ be the subgroup formed by the principal fractional ideals whose generators reduce to the identity modulo $\m$. 
By definition, the  {\em ray class group modulo $\m$} is $\Cl_\m : = \I_\m /i(K_{\m,1})$. 
Taking  $\Gamma$ into account now we write $i(K_{\m,\Gamma})$ for the subgroup of $\I_\m$ consisting of principal fractional ideals
whose generators reduce to an element of $\Gamma$ modulo $\m$. The  quotient 
$\I_\m / i(K_{\m,\Gamma}) \cong  \Cl_\m  /  (i(K_{\m,\Gamma})/ i(K_{\m,1}))$, which  can be regarded as a generalized `class group',  plays a key role in the parametrization of  low-temperature KMS states in \cite[Theorem~3.2(iii)]{Bru2} and in the K-theory formula in \cite[Theorem 4.1]{BruLi}.
 
The congruence monoid  $R_{\m,\Gamma}$ acts on $R$ by multiplication and gives rise to the semi-direct product $R\rtimes R_{\m,\Gamma}$, which is a submonoid of the full $ax+b$-semigroup $\rxrx$. 
By \cite[Proposition 3.3]{Bru1}, $R\rtimes R_{\m,\Gamma}$ is a left Ore monoid that has $Q\rtimes K_{\m,\Gamma}$ as its group of left quotients, where $Q:=R_\m^{-1}R$ denotes the localization of $R$ at $R_\m$ as usual.

We are interested in the C*-dynamical system consisting of the left regular C*-algebra $C_\lambda^*(R\rtimes R_{\m,\Gamma})$ 
under the canonical time evolution $\sigma$ given by
$\sigma_t(\lambda_{(b,a)})=N(a)^{it}\lambda_{(b,a)}$ for all $(b,a)\in R\rtimes R_{\m,\Gamma}$, see \cite[Proposition~3.1]{Bru2}.  
Fix an integral ideal $\a_\k\in\k$ for each class $\k\in \I_\m /i(K_{\m,\Gamma})$ (for the trivial class take $\a_{[R]}=R$). According to \cite[Theorem~3.2(iii)]{Bru2}, for each $\beta \in (2,\infty)$ there is
an affine isomorphism of the simplex of tracial states on the C*-algebra 
\[
\bigoplus_{\k\in \I_\m / i(K_{\m,\Gamma})}C^*(\a_\k\rtimes R_{\m,\Gamma}^*)
\]
and the simplex of KMS$_\beta$ states on $C_\lambda^*(R\rtimes R_{\m,\Gamma})$.

Recall that $C_\lambda^*(R\rtimes R_{\m,\Gamma})$ is canonically isomorphic to the C*-algebra of
the partial transformation groupoid $(Q\rtimes K_{\m,\Gamma})\ltimes\Omega$ 
given in equation (3) in \cite[\S~5.2]{Bru1}.
 The unit space $\Omega:=\Omega_R^\m$ of the groupoid is best viewed as an adelic space on which
  $Q\rtimes K_{\m,\Gamma}$ acts canonically, see  \cite[\S~5.2]{Bru1}. \tcb{This adelic space is defined as follows. For each non-zero prime ideal $\p$ of $R$, let $K_\p$ denote the $\p$-adic completion of $K$, and let $R_\p$ be the ring of $\p$-adic integers in $K_\p$. Then the compact ring $\hat{R}_S:=\prod_{\p\nmid \m_0}R_\p$ has unit group $\hat{R}_S^*=\prod_{\p\nmid \m_0}R_\p^*$, and we let $\bar{\mathbf{a}}$ denote the image of $\mathbf{a}\in\hat{R}_S$ under the quotient map $\hat{R}_S\to \hat{R}_S/\hat{R}_S^*$. Then $\Omega$ is defined as the quotient of $\hat{R}_S\times \hat{R}_S/\hat{R}_S^*$ by the equivalence relation $(\mathbf{b},\bar{\mathbf{a}})\sim (\mathbf{d},\bar{\mathbf{c}})$ if and only if $\bar{\mathbf{a}}=\bar{\mathbf{c}}$ and $\mathbf{b}-\mathbf{d}\in \mathbf{a}\hat{R}_S$. We let $[\mathbf{b},\bar{\mathbf{a}}]$ denote the equivalence class of $(\mathbf{b},\bar{\mathbf{a}})\in \hat{R}_S\times \hat{R}_S/\hat{R}_S^*$. The partial action of $Q\rtimes K_{\m,\Gamma}$ on $\Omega$ then comes from the canonical partial action of $Q\rtimes K_{\m,\Gamma}$ on $\hat{R}_S\times \hat{R}_S/\hat{R}_S^*$.}
  
Since for $\beta  > 2$ the extremal quasi-invariant measures \tcb{on $\Omega$} are supported on orbits, \cite[Corollary~1.4]{Nesh} applies and prescribes that KMS$_\beta$ states are parametrized by the tracial states on the isotropy group of a point in each orbit.
In order to construct the state $\varphi_{\beta,\k,\tau}$ associated to a trace $\tau$ of $C^*(\a_\k\rtimes R_{\m,\Gamma}^*)$ 
 one needs to find a quasi-invariant measure supported on the orbit of $[0,\a_\k] \in \Omega$, show that the corresponding isotropy subgroup is $\a_\k\rtimes R_{\m,\Gamma}^*$ and then transfer the trace $\tau$ of $C^*(\a_\k\rtimes R_{\m,\Gamma}^*)$ to the isotropy subgroups of the rest of the orbit, which are obtained via conjugation by elements of the enveloping group
 $Q \rtimes K_{\m,\Gamma}$. We set the notation and give the precise statements in the following propositions.

\begin{proposition}\label{pro:orbitandisotropy}
Let $(Q\rtimes K_{\m,\Gamma})\ltimes\Omega$ be the partial transformation groupoid from \cite[\S~5.2]{Bru1}. 
 As before, fix an integral ideal $\a_\k\in\k$ for each class $\k\in \I_\m /i(K_{\m,\Gamma})$ (for the trivial class take $\a_{[R]}=R$)
 and for each $\a \in \k$ fix a generator $\sofa \in K_{\m,\Gamma}$ of the principal ideal $\a\inv\a_\k$ (for $\k=[R]$ and $\a=R$, take $\sofa =1$).
Using conjugation in $Q \rtimes K_{\m,\Gamma}$ as maps between its subgroups, we then have
\begin{enumerate}
\item[(i)] the orbit of $[0,\a_\k]$ is $\{ [x, \sofa\inv \a_\k] \in\Omega: x\in R,  \ \a \in \k\}$, and the isotropy at $[x, \sofa\inv \a_\k]$ is the group   
$G_{x,\a}:=(x,1)(\sofa\inv \a_\k\rtimes R_{\m,\Gamma}^*)(-x,1)$;

\item[(ii)] the isomorphisms of C*-algebras of isotropy groups derived from the partial action of $Q \rtimes K_{\m,\Gamma}$ are given by
\[
\theta_{x,\a}: u_{(r,g)} \mapsto u_{(t_\a (r+(g-1)x), g)}
\qquad \quad C^*(G_{x,\a}) \cong C^*(\a_\k\rtimes R_{\m,\Gamma}^*),
\]
where  $u_{(r,g)}$ denotes the canonical unitary generator of $C^*(G_{x,\a})$
corresponding to $(r,g) \in G_{x,\a}$;

\item[(iii)]  the corresponding simplices of tracial states are affinely isomorphic through the maps
\[
\tau \mapsto \tau_{x,\a} := \tau\circ \theta_{x,\a} \qquad \quad T(C^*(\a_\k\rtimes R_{\m,\Gamma}^*)) \to T(C^*(G_{x,\a})),
\]
where for each $ \tau \in T(C^*(\a_\k\rtimes R_{\m,\Gamma}^*))$ the trace $\tau_{x,\a}$ is characterized by  
\[
\tau_{x,\a}(u_{(x,1)}u_{(r,g)}u_{(x,1)}^*)=\tau(u_{(\sofa r,g)}) \qquad (r,g)\in \a \rtimes R_{\m,\Gamma}^*.
\]
\end{enumerate}
\end{proposition}
\begin{proof}
The proof consists of several straightforward calculations which we omit. 
\end{proof}

For each class $\k\in \I_\m /i(K_{\m,\Gamma})$, let 
\[
\zeta_{K,\k}(s):=\sum_{\a\in\k,\a\subseteq R}N(\a)^{-s}
\]
denote the \emph{partial zeta function associated with $\k$}. The series $\zeta_{K,\k}(s)$ converges for all $s\in\CC$ with $\Re s>1$, see, for example, \cite[Ch.VI Section 2]{MilCFT}. We will write $\zeta_\k(s)$ rather than $\zeta_{K,\k}(s)$ when it will not cause confusion.

 Using the quasi-invariant probability measures supported on orbits determined in \cite[\S~3]{Bru2} and following the construction in the proof of \cite[Theorem 1.3]{Nesh} we  derive an explicit formula for the low-temperature KMS states as follows.
 
 \begin{proposition}
 \label{pro:KMSfromtrace}
Fix $\beta>2$. If $\tau$ is a tracial state on $C^*(\a_\k\rtimes R_{\m,\Gamma}^*)$, then the KMS$_\beta$ state $\varphi_{\beta,\k,\tau}$ is given by
\begin{equation}
\label{eqn:KMSformula1}
\varphi_{\beta,\k,\tau}(f)=\frac{1}{\zeta_\k(\beta-1)}\sum_{\a\in\k, \a\subseteq R}\sum_{x\in R/\a}\sum_{\gamma\in G_{x,\a}}N(\a)^{-\beta}f((\gamma,[x,\a]))\tau_{x,\a}(u_\gamma)
\end{equation}
for all $f\in C_c((Q\rtimes K_{\m,\Gamma})\ltimes\Omega)$.
 \end{proposition}
 \begin{proof}
From \cite[Proof~of~Theorem~3.2(iii)]{Bru2}, we see that the quasi-invariant probability measure on $\Omega$ determined by the restriction of $\varphi_{\beta,\k,\tau}$ to $C(\Omega)$ is given by $\mu_{\beta,\k}=\pi_*(m\times\nu_{\beta-1,\k})$ where $m$ is the normalized Haar measure on $\hat{R}_S$, $\nu_{\beta-1,\k}$ is the probability measure on $\hat{R}_S/\hat{R}_S^*\cong\prod_{\p\nmid \m_0}\p^{\NN\cup\{\infty\}}$ given by $\nu_{\beta-1,\k}=\frac{1}{\zeta_\k(\beta-1)}\sum_{\b\in\k, \b\subseteq R}N(\b)^{-(\beta-1)}\delta_\b$, and $\pi:\hat{R}_S\times \hat{R}_S/\hat{R}_S^*\to\Omega$ is the quotient map $\pi(\mathbf{b},\mathbf{a})=[\mathbf{b},\bar{\mathbf{a}}]$.

We now show that $\mu_{\beta,\k}$ has the following explicit description:
\begin{equation}
\label{eqn:measures}
\mu_{\beta,\k}=\frac{1}{\zeta_\k(\beta-1)}\sum_{\b\in\k,\b\subseteq R}N(\b)^{-\beta}\sum_{y\in R/\b}\delta_{[y,\b]}.
\end{equation}
For this, it suffices to show that these measures agree on the following compact open subsets of $\Omega$ given by 
\[
V_{x+\a}:=\{[\mathbf{b},\bar{\mathbf{a}}] : \a\mid \bar{\mathbf{a}} \text{ and } \a\mid (\mathbf{b}-x)\} \quad \text{for $x\in R$ and integral $\a\in\I_\m$}
\]
because the characteristic functions $\{1_{V_{x+\a}} : x\in R,\a\in\I_\m,\a\subseteq R\}$ span a dense *-subalgebra of $C(\Omega)$.
We have
\begin{align*}
\mu_{\beta,\k}(V_{x+\a})&=m\times\nu_{\beta-1,\k}(\pi^{-1}V_{x+\a}) \\
&=m\times\nu_{\beta-1,\k}((x+\a\hat{R}_S)\times \a\hat{R}_S/\hat{R}_S^*)\quad\text{ since $\pi^{-1}(V_{x+\a})=(x+\a\hat{R}_S)\times \a\hat{R}_S/\hat{R}_S^*$}\\
&=m(\a\hat{R}_S)\nu_{\beta-1,\k}(\a\hat{R}_S/\hat{R}_S^*)\\
&=\frac{N(\a)^{-1}}{\zeta_\k(\beta-1)}\sum_{\b\in\k, \a\mid \b}N(\b)^{-(\beta-1)}.
\end{align*}
For an integral ideal $\b$ in $\I_\m$ and $y\in R/\b$, we have $[y,\b]\in V_{x+\a}$ if and only if $\a\mid\b$ and $y-x\in\a$, so 
\[
\delta_{[y,\b]}(V_{x+\a})=\begin{cases}
1 & \text{ if } \a\mid \b \text{ and } y-x\in \a,\\
0 & \text{ otherwise.}
\end{cases}
\]
If $\a\mid \b$, then there are $|\a/\b|=N(\b)/N(\a)$ elements $y\in R/\b$ such that $y-x\in\a$, so we have
\[
\frac{1}{\zeta_\k(\beta-1)}\sum_{\b\in\k,\b\subseteq R}N(\b)^{-\beta}\sum_{y\in R/\b}\delta_{[y,\b]}(V_{x+\a})=\frac{1}{\zeta_\k(\beta-1)}\sum_{\b\in\k,\a\mid\b}N(\b)^{-\beta}N(\b)N(\a)^{-1},
\]
which shows that \eqref{eqn:measures} holds.

Since $\mu_{\beta,\k}$ is concentrated on the orbit $\{[x,\a] : \a\in\k, \a\subseteq R, x\in R\}$ for the partial action $Q\rtimes K_{\m,\Gamma}\curvearrowright\Omega$, $\tau\in T(C^*(\a_\k\rtimes R_{\m,\Gamma}^*))$ canonically gives rise to a $\mu_{\beta,\k}$-measurable field of tracial states (see \cite[Proof~of~Corollary~1.4]{Nesh}); namely, the $\mu_{\beta,\k}$-measurable field corresponding to $\tau$ is $(\mu,\{\tau_{x,\a}\})$ where we use that condition (iii) in Proposition~\ref{pro:orbitandisotropy} is exactly condition (ii) of \cite[Theorem 1.3]{Nesh}.
 The KMS state $\varphi_{\beta,\k,\tau}$ is then the KMS state corresponding to  $(\mu,\{\tau_{x,\a}\})$ via \cite[Theorem~1.3]{Nesh}, which is given explicitly by
\[
\varphi_{\beta,\k,\tau}(f)=\sum_{\substack{[x,\a]\in\Omega\\  \a\in\k,\a\subseteq R, x\in R}}N(\a)^{-\beta}\sum_{\gamma\in G_{x,a}}f((\gamma,[x,\a]))\tau_{x,a}(u_\gamma).
\]
for all $f\in C_c((Q\rtimes K_{\m,\Gamma})\ltimes\Omega)$. Note that the explicit formula is not given in the statement of \cite[Theorem~1.3]{Nesh}, so one must look into the proof (see also \cite[Theorem~1.1]{LLN2} where a version of \cite[Theorem~1.3]{Nesh} with the explicit formula is stated). Since this formula coincides with the one given in \eqref{eqn:KMSformula1}, we are done.
 \end{proof}

Next we characterize these low-temperature KMS$_\beta$ states $\varphi_{\beta,\k,\tau}$ by their values on a set of analytic elements spanning a dense *-subalgebra of $C^*_\lambda(R \rtimes R_{\m,\Gamma})$ with a formula similar to the one given in \cite[Remark~7.4(2)]{CDL}. This will be instrumental later when we construct suitable Hilbert spaces for the computation of the partition functions of the various states following Proposition~\ref{prop:gibbsbydensity}. 

To lighten the notation we will use the following conventions. For $b\in R_{\m,\Gamma}$, we write $s_b:=\lambda_{(0,b)}$; for $x\in R$ we write  $u^x:=\lambda_{(x,1)}$, so that $\lambda_{(x,b)} = u^x s_b$, and for a non-empty constructible right ideal $(y+\b)\times (\b\cap R_{\m,\Gamma})$ of $R\rtimes R_{\m,\Gamma}$, we write $e_{y+\b}:=E_{(y+\b)\times (\b\cap R_{\m,\Gamma})}$ (see \cite[\S~3~\&~\S~4]{Bru1}).

\begin{proposition}
\label{pro:KMSformula}
The linear span of the elements  $s_b^*e_{y+\b}u^ds_c$ with  $b,c\in R_{\m,\Gamma}$, $y,d\in R$,  and $\b\in\I_\m$ with $\b\subseteq R$ is a dense *-algebra of $C_\lambda^*(R\rtimes R_{\m,\Gamma})$. Suppose $\beta \in (2,\infty)$ and let $\tau $  be a tracial state of $C^*(\a_\k\rtimes R_{\m,\Gamma}^*)$ for a class $\k \in \I_\m/i(K_{\m,\Gamma})$.  Then 
\begin{equation}\label{eq:kmsvalues}
\varphi_{\beta,\k,\tau}(s_b^*e_{y+\b}u^ds_c)=\begin{cases}\displaystyle
\frac{1}{\zeta_\k(\beta-1)}\sum_{\substack{\a\in\k,\\ b\a\subseteq \b}}\ \sum_{\substack{x\in R/\a,\; d+cx-y\in\b,\\ d+bx(g-1)\in b\a}}N(\a)^{-\beta}\tau(u_{(t_\a(b^{-1}d+x(g-1)),g)}) &\text{ if } g:=b^{-1}c\in R_{\m,\Gamma}^*  \\
\quad 0 &\text{ otherwise}
\end{cases}
\end{equation}
\end{proposition}
\begin{proof}
Since $\lambda_{(x,b)}^* e_{y+ \b} \lambda_{(d,c)} = s_b^* u^{-x} e_{y+ \b} u^d s_c = s_b^* e_{ y-x + \b} u^{d-x} s_c$ 
 a computation similar to the one given in \cite[Remark 1.6]{La00} for general Ore semigroups shows that
 the set of products of the form $s_b^*e_{y+\b}u^ds_c$ is closed under multiplication.
Thus their linear span is a *-subalgebra that  contains the generating elements $s_b$,  $u^x$, and $e_\b$, hence is dense in 
$C_\lambda^*(R\rtimes R_{\m,\Gamma})$.

 To verify formula \eqref{eq:kmsvalues} we use the isomorphism of $C_\lambda^*(R\rtimes R_{\m,\Gamma})$ to a partial transformation groupoid C*-algebra, in which  the spanning element becomes a convolution of functions of compact support: 
 \[
 s_b^*e_{y+\b}u^ds_c = 1_{\{(0,b)\}\times\Omega}^**\hat{e}_{y+\b}*1_{\{(d,c)\}\times\Omega},
 \]
 where $\hat{e}_{y+\b}\in C(\Omega)$ denotes the image of the projection $e_{y+\b}\in C_\lambda^*(R\rtimes R_{\m,\Gamma})$ under the isomorphism $C_\lambda^*(R\rtimes R_{\m,\Gamma})\cong C_r^*((Q\rtimes K_{\m,\Gamma})\ltimes\Omega)$.
 
Let $(\gamma,[\textbf{b},\bar{\textbf{a}}])\in (Q\rtimes K_{\m,\Gamma})\ltimes\Omega$ with $\gamma=(n,k)$ and $[\textbf{b},\bar{\textbf{a}}]\in \Omega$. Then
\begin{align*}
&(1_{\{(0,b)\}\times\Omega}^**\hat{e}_{y+\b}*1_{\{(d,c)\}\times\Omega})(\gamma,[\textbf{b},\bar{\textbf{a}}])\\
&=\sum_{(h,[\textbf{d},\bar{\textbf{c}}])\in\G^{\gamma.[\textbf{b},\bar{\textbf{a}}]}}1_{\{(0,b)\}\times\Omega}((h^{-1},h[\textbf{d},\bar{\textbf{c}}]))(\hat{e}_{y+\b}*1_{\{(d,c)\}\times\Omega})((h^{-1},h[\textbf{d},\bar{\textbf{c}}])(\gamma,[\textbf{b},\bar{\textbf{a}}]))\\
&=(\hat{e}_{y+\b}*1_{\{(d,c)\}\times\Omega})((0,b)\gamma,[\textbf{b},\bar{\textbf{a}}])\\
&=\hat{e}_{y+\b}([bn+bk\textbf{b},bk\bar{\textbf{a}}])1_{\{(d,c)\}\times\Omega}((bn,bk),[\textbf{b},\bar{\textbf{a}}])\\
&=\begin{cases}
1 & \text{ if } d=bn,\; c=bk,\; \b\mid bk\bar{\textbf{a}},\;\text{ and } \b\mid (bn+bk\textbf{b}-y), \\
0 & \text{ otherwise,}
\end{cases}\\
&=\begin{cases}
1 & \text{ if } (b^{-1}d,b^{-1}c)=(n,k),\; \b\mid c\bar{\textbf{a}},\;\text{ and } \b\mid (d+c\textbf{b}-y), \\
0 & \text{ otherwise.}
\end{cases}
\end{align*}

Suppose now that $\a$ is an integral ideal in the class $\k$, let $x\in R$ so that $[\textbf{b},\bar{\textbf{a}}] = [x,\a]$
is a point in the orbit of $[0,\a_\k]$. Assuming $\gamma=(n,k)$ is in the corresponding isotropy subgroup  $G_{x,\a}=(x,1)\a\rtimes R_{\m,\Gamma}^*(-x,1)$ we may write it as $\gamma=(x,1)(r,g)(-x,1)=(r+x(1-g),g)$ for some $(r,g)\in\a\rtimes R_{\m,\Gamma}^*$.  Then
\begin{align*}
(1_{\{(0,b)\}\times\Omega}^**\hat{e}_{y+\b}*&1_{\{(d,c)\}\times\Omega})(\gamma,[x,\a])\\
&=\begin{cases}
1 & \text{ if } (b^{-1}d,b^{-1}c)=(r+x(1-g),g),\; c\a\subseteq \b,\;\text{ and } d+cx-y\in\b, \\
0 & \text{ otherwise,}
\end{cases}\\
&=\begin{cases}
1 & \text{ if } d+bx(g-1)=br,\; b^{-1}c=g,\; b\a\subseteq \b,\;\text{ and } d+cx-y\in\b, \\
0 & \text{ otherwise.}
\end{cases}
\end{align*}
If $b^{-1}c$ is not in $R_{\m,\Gamma}^*$, the function $1_{\{(0,b)\}\times\Omega}^**\hat{e}_{y+\b}*1_{\{(d,c)\}\times\Omega}$ vanishes on the orbit of $[0,\a_\k]$ so the right hand side of equation~\eqref{eqn:KMSformula1} vanishes.
If $g = b^{-1}c \in R_{\m,\Gamma}^*$, then equation~\eqref{eqn:KMSformula1} yields
\begin{align*}
\varphi_{\beta,\k,\tau}(s_b^*e_{y+\b}u^ds_c)&=\frac{1}{\zeta_\k(\beta-1)}\sum_{\substack{\a\in\k,\\ b\a\subseteq \b}}\ \sum_{\substack{x\in R/\a,\; d+cx-y\in\b,\\ d+bx(g-1)\in b\a}}N(\a)^{-\beta}\tau_{x,\a}(u_{(x,1)}u_{(b^{-1}d+x(g-1),g)}u_{(x,1)}^*)
\\
&=\frac{1}{\zeta_\k(\beta-1)}\sum_{\substack{\a\in\k,\\ b\a\subseteq \b}}\ \sum_{\substack{x\in R/\a,\; d+cx-y\in\b,\\ d+bx(g-1)\in b\a}}N(\a)^{-\beta}\tau(u_{(t_\a(b^{-1}d+x(g-1)),g)}),
\end{align*}
where in the last step we have used Proposition~\ref{pro:orbitandisotropy}(iii). This completes the proof of formula \eqref{eq:kmsvalues}.
\end{proof}

\section{Tracial states on $C^*(\a_\k\rtimes R_{\m,\Gamma}^*)$}
\label{sec:traces}
Given a C*-algebra $A$, we equip  the set  $T(A)$ of tracial states on $A$  with the weak-* topology. Given a convex set $\Sigma$, we denote by $\ex\Sigma$ the set of extreme points in $\Sigma$.
We will need  the following general lemma about traces on direct sums, which we suspect to be well-known, but we provide the proof for completeness.
\begin{lemma}
\label{lem:tracesonsum}
Let $A_1,A_2,...,A_n$ be C*-algebras, and let $A:=\bigoplus_{k=1}^nA_k$. Then each  projection map $\pi_k:A\to A_k$ induces a continuous injection of extremal tracial states $\pi_k^*:\ex T(A_k)\to \ex T(A)$ such that $\pi_k^*(\tau)(a)=\tau(\pi_k(a))$, and these maps induce a homeomorphism
\[
\sqcup_{k=1}^n \pi_k^*: \bigsqcup_{k=1}^n\ex T(A_k) \simeq \ex T(A).
\]
\end{lemma}
\begin{proof}
It is easy to see that each $\pi_k^*$ is a continuous injection, so that $\sqcup_{k=1}^n \pi_k^*: \bigsqcup_{k=1}^n\ex T(A_k)\to \ex T(A)$ is also a continuous injection. It is not difficult to see that $\sqcup_{k=1}^n \pi_k^*$ is also surjective. To show that $\sqcup_{k=1}^n \pi_k^*$ is a homeomorphism, it is enough to show that $\pi_k^*(\ex T(A_K))$ is open in $\ex T(A)$ for each $k$.

For $\tau\in\ex T(A)$, the restriction $\tau\vert_{A_k}$ of $\tau$ to $A_k$ either lies in $\ex T(A_k)$ or is the zero functional, so we get a continuous map
\[
\res_A^{A_k}: \ex T(A)\to \ex T(A_k)\cup\{0\}.
\]
Hence, the set 
\[
(\res_A^{A_k})^{-1}(0)=\bigsqcup_{l\neq k}\pi_l^*(\ex T(A_l))
\]
is closed, which implies that $\pi_k^*(\ex T(A_k))$ is open.
\end{proof}

In light of Lemma~\ref{lem:tracesonsum}, we see that if a trace $\tau\in T(C^*(\a_\k\rtimes R_{\m,\Gamma}^*))$ is extremal, then the KMS$_\beta$ state $\varphi_{\beta,\k,\tau}$ from Proposition~\ref{pro:KMSformula} is extremal and is thus a factor state.
In order to determine the type of this factor we will need to obtain a more explicit description of the tracial states of the group C*-algebras $C^*(\a_\k\rtimes R_{\m,\Gamma}^*)$ for $\k\in \I_\m/i(K_{\m,\Gamma})$.

The proofs of \cite[Lemma~2.1]{LaWa} and \cite[Theorem~2.2]{LaWa},  with the group of units replaced by its finite index subgroup $R_{\m,\Gamma}^*$, carry over to the present setting, showing that extremal traces arise from ergodic measures of groups of linear toral automorphisms and characters of the associated isotropy subgroups. 
\begin{theorem}
\label{thm:tracesandmeasures}
Let $\k \in \I_\m/i(K_{\m,\Gamma})$ and let $R_{\m,\Gamma}^*\curvearrowright \widehat\a_\k$ be the action of 
$R_{\m,\Gamma}^*$ on the Pontryagin dual $\widehat\a_\k$ of the integral ideal $\a_\k$ obtained via transposition by duality from the multiplicative action of $R_{\m,\Gamma}^*$ on  $\a_\k$.
\begin{enumerate}
\item[(i)] For any ergodic $R_{\m,\Gamma}^*$-invariant probability measure $\mu$ on $\widehat{\a_\k}$, the system $R_{\m,\Gamma}^*\curvearrowright (\widehat{\a_\k},\mu)$ has constant isotropy $\mu$-a.e., that is, there exists a unique subgroup $H_\mu$ of $R_{\m,\Gamma}^*$ such that $H_\mu=(R_{\m,\Gamma}^*)_\chi=\{v\in R_{\m,\Gamma}^* : v.\chi=\chi\}$ for $\mu$-a.e. $\chi\in\widehat{\a_\k}$.
\item[(ii)] Let $u_{(b,v)}$ denote the canonical unitary in $C^*(\a_\k\rtimes R_{\m,\Gamma}^*)$ corresponding to $(b,v)\in \a_\k\rtimes R_{\m,\Gamma}^*$. For each extremal tracial state $\tau$ of $C^*(\a_\k\rtimes R_{\m,\Gamma}^*)$, there exists a unique probability \tcb{measure} $\mu_\tau$ on $\widehat{\a_\k}$ such that
\begin{equation}\label{eq:tracefrommeasure}
\tau(u_{\tcb{(b,1)}})=\int_{\widehat{\a_\k}} \gamma(b)d\mu_\tau(\gamma) \quad \text{for } b\in\a_\k. 
\end{equation}
The measure $\mu_\tau$ is ergodic and $R_{\m,\Gamma}^*$-invariant, and the function $\chi_\tau$ defined by $\chi_\tau(h)=\tau(u_h)$ for $h\in H_{\mu_\tau}$ is a character on $H_{\tcb{\mu_\tau}}$.
\item[(iii)] The map $\tau\mapsto (\mu_\tau,\chi_\tau)$ is a bijection from the set of extremal tracial states of $C^*(\a_\k\rtimes R_{\m,\Gamma}^*)$ onto the set of pairs $(\mu,\chi)$ consisting of an ergodic $R_{\m,\Gamma}^*$-invariant probability measure $\mu$ on $\widehat{\a_\k}$ and a character $\chi\in\widehat{H_\mu}$.
\end{enumerate} 
\end{theorem}

We shall also need to know the type of the factors associated to extremal tracial states on each of the direct summands of $\bigoplus_\k C^*(\a_\k\rtimes R_{\m,\Gamma}^*)$.

\begin{proposition}
\label{prop:typeoftraces}
Let $\tau$ be an extremal tracial state on $C^*(\a_\k\rtimes R_{\m,\Gamma}^*)\cong C(\widehat{\a_\k})\rtimes R_{\m,\Gamma}^*$ and let $(\mu,\chi)=(\mu_\tau,\chi_\tau)$ be the pair corresponding to $\tau$ under Theorem~\ref{thm:tracesandmeasures}.   If the measure $\mu$ is concentrated on a finite orbit $O$, then $\pi_\tau(C^*(\a_\k\rtimes R_{\m,\Gamma}^*))''$ is a factor of type I$_{|O|}$, otherwise it is
a factor of type II$_1$.

For each pair  $(O,\chi)$ consisting of  a
finite orbit $O$ for $R_{\m,\Gamma}^*\curvearrowright \widehat{\a_\k}$ and a character $\chi$ of  the isotropy group $(R_{\m,\Gamma}^*)_O$ of any point of $O$,  define
\[
\tau_{O,\chi}(fu_g)=\begin{cases}
\frac{\chi(g)}{|O|}\sum_{x\in O}f(x) & \text{ if } g\in (R_{\m,\Gamma}^*)_O\\
0 & \text{ otherwise}
\end{cases}
\]
for $g\in R_{\m,\Gamma}^*$ and $f\in C(\widehat{\a_\k})$.
Then the map $(O,\chi)\mapsto \tau_{O,\chi}$ is a bijection between the set of such pairs and  the set of extremal tracial states of type I on $C^*(\a_\k\rtimes R_{\m,\Gamma}^*)$. 
\end{proposition}
\begin{proof}
Let $\tilde{\chi}$ be an extension of $\chi$ to $R_{\m,\Gamma}^*$. Then there is an isomorphism $\pi_\tau(C^*(\a_\k\rtimes R_{\m,\Gamma}^*))''\cong L^\infty(\widehat{\a}_\k,\mu)\rtimes (R_{\m,\Gamma}^*/H_\mu)$ such that for $f\in C(\widehat{\a}_\k)$, $\pi_\tau(f)\mapsto f\in L^\infty(\widehat{\a}_\k,\mu)$ and $\pi_\tau(u_w)\mapsto \tilde{\chi}(w)u_{\bar{w}}$ where $\bar{w}$ is the image of $w$ under the quotient map $R_{\m,\Gamma}^*\to R_{\m,\Gamma}^*/H_\mu$ (see \cite[Remark~2.5]{Nesh}). The type assertions now follow from the Murray--von Neumann  classification of transformation group von Neumann algebras \cite[Lemma 13.1.1 and Lemma 13.1.2]{MvN1}.

The equi-probability on the points of a finite orbit is the only ergodic invariant measure supported on the orbit, 
so the second claim follows from Theorem~\ref{thm:tracesandmeasures}, 
and equation \eqref{eq:tracefrommeasure} gives the formula for $\tau_{O,\chi}$.
\end{proof}

\begin{remark}
\label{rmk:Furstenberg}
When $K=\QQ$ or $K$ is quadratic imaginary the group  $R^*$ of units is finite
and the same holds for $R_{\m,\Gamma}^*$, hence all the ergodic invariant measures
of $R_{\m,\Gamma}^*\curvearrowright \widehat\a_\k$
are supported on finite orbits and all extremal traces are of type I. If the rank of $R^*$ is at least 1, then so is that of $R^*_{\m,\Gamma}$ and hence normalized Haar measure on $\widehat\a_\k$ is ergodic invariant. 
Whether there are other ergodic invariant probability measures than those arising from finite orbits and from Haar measure is
the generalized Furstenberg question, which is open for non CM fields with unit rank 2 or higher \cite{LaWa}.
In any case, for all fields, the finite orbits exhaust type I extremal traces, so Haar measure and the possible extra ones are
necessarily of type II$_1$.
\end{remark}

\section{Type classification for low-temperature KMS states}
\label{sec:type}
The two main goals of this section are to determine the type of the extremal KMS$_\beta$ states for $\beta\in (2,\infty)$,  Corollary~\ref{cor:typeKMS}, and to compute the partition functions for those of type I,  Theorem~\ref{thm:partitionfncI}. 
Both goals depend heavily on  Proposition~\ref{prop:extendrep}, in which we `induce' representations of the isotropy group C*-algebras $C^*(\a_\k\rtimes R_{\m,\Gamma}^*)$ to representations of $C_\lambda^*(R\rtimes R_{\m,\Gamma})$. The type of the `induced' representation can then be  computed in terms of that of the original representation,  Theorem~\ref{thm:vNalg}. Moreover, when we induce from a type I extremal tracial state
of $C^*(\a_\k\rtimes R_{\m,\Gamma}^*)$, the construction gives rise to an irreducible representation 
and thus to an admissible triple, which allows us to compute the partition function of the associated KMS$_\beta$ state. This construction is inspired by \cite[Lemma~7.1]{CDL}. Our notation for representations induced via Hilbert bimodules is the one used in \cite{RaeWill}.

As before we use $Q:=R_\m^{-1}R$ to denote the localization of $R$ at $R_\m$, so that 
$Q\rtimes K_{\m,\Gamma}$ is the group of left quotients of $R\rtimes R_{\m,\Gamma}$ by \cite[Proposition~3.3]{Bru1}.
  Denote by $u_\gamma$ the canonical generating unitary in $C^*(Q\rtimes K_{\m,\Gamma})$ corresponding to $\gamma\in Q\rtimes K_{\m,\Gamma}$.
For each fractional ideal $\a\in\I_\m$, the group $\a\rtimes R_{\m,\Gamma}^*$ is a subgroup of $Q\rtimes K_{\m,\Gamma}$, and thus  the C*-algebra $C^*(\a\rtimes R_{\m,\Gamma}^*)$ can be viewed as the subalgebra $C^*(\{u_\gamma: \gamma \in \a\rtimes R_{\m,\Gamma}^*\})$ of $C^*(Q\rtimes K_{\m,\Gamma})$.

Recall from Proposition~\ref{pro:orbitandisotropy} that we have chosen a reference ideal $\a_\k$ in each class $\k \in \I_\m /i(K_{\m,\Gamma})$, with the trivial class being represented by $[R]$ itself,  and for each $\a \in \k$ we have also chosen $\sofa \in K_{\m,\Gamma}$ such that 
$\sofa \a = \a_\k$.  These choices give specific isomorphisms $\theta_{\sofa} : C^*(\a\rtimes R_{\m,\Gamma}^*) \to C^*(\a_\k\rtimes R_{\m,\Gamma}^*)$ between subalgebras of $C^*(Q\rtimes K_{\m,\Gamma})$ determined by  $\theta_{\sofa}(u_{(x,v)})=u_{(0,\sofa)}u_{(x,v)}u_{(0,\sofa)}^*=u_{(\sofa x,v)}$ for all $(x,v)\in\a\rtimes R_{\m,\Gamma}^*$.

Let $\k\in \I_\m/i(K_{\m,\Gamma})$. The canonical right Hilbert $C^*(\a_\k\rtimes R_{\m,\Gamma}^*)$-module structure on $C^*(Q\rtimes K_{\m,\Gamma})$ arising from the inclusion of groups $\a_\k\rtimes R_{\m,\Gamma}^*\subseteq Q\rtimes K_{\m,\Gamma}$ gives rise to a right Hilbert $C^*(\a_\k\rtimes R_{\m,\Gamma}^*)$-module structure on 
\[
\X_\a:=C^*(R\rtimes R_{\m,\Gamma}^*)u_{(0,\sofa)}^*\subseteq C^*(Q\rtimes K_{\m,\Gamma})
\] 
with $C^*(\a_\k\rtimes R_{\m,\Gamma}^*)$-valued inner product determined by 
\[
\langle u_{(y,h)}u_{(0,t_\a)}^*,u_{(x,g)}u_{(0,t_\a)}^*\rangle=\begin{cases}
u_{(0,t_\a)}u_{(y,h)}^*u_{(x,g)}u_{(0,t_\a)}^* & \text{ if } (y,h)^{-1}(x,g)\in \a\rtimes R_{\m,\Gamma}^*,\\
0 & \text{ if }(y,h)^{-1}(x,g)\notin \a\rtimes R_{\m,\Gamma}^*.
\end{cases}
\]

We shall view $\X_\a$ as a $C^*(R\rtimes R_{\m,\Gamma}^*)$--$C^*(\a_\k\rtimes R_{\m,\Gamma}^*)$-bimodule via the canonical left action of $C^*(R\rtimes R_{\m,\Gamma}^*)$ given by multiplication.

We introduce now some notation that will be used throughout the remainder of this section. For each integral ideal $\a\in\k$, let $\R_\a$ be a complete set of representatives for $R/\a$ such that $0\in\R_\a$. For $x\in R$, let $\bar{x}^\a$ denote the unique element of $\R_\a$ such that $x-\bar{x}^\a\in\a$. We will occasionally omit the superscript when the notation becomes particularly cumbersome.

Let  $\k$ be a class in $\I_\m/i(K_{\m,\Gamma})$. 
Suppose $\pi_0$ is a representation of $C^*(\a_\k\rtimes R_{\m,\Gamma}^*)$ on a Hilbert space $\H_0$, and for each integral ideal $\a\in\k$ let 
\[
\vartheta_\a:=\X_\a \dash\Ind_{C^*(\a_\k\rtimes R_{\m,\Gamma}^*)}^{C^*(R\rtimes R_{\m,\Gamma}^*)}\pi_0
\]
 be the induced representation of $C^*(R\rtimes R_{\m,\Gamma}^*)$ on the Hilbert space 
\[
\H_\a:=\X_\a\otimes_{C^*(\a_\k\rtimes R_{\m,\Gamma}^*)} \H_0\subseteq C^*(Q\rtimes K_{\m,\Gamma})\otimes_{C^*(\a_\k\rtimes R_{\m,\Gamma}^*)} \H_0.
\]
\begin{lemma}
\label{lem:inducedrep} Let $\a$ be an integral ideal in $\I_\m$, let $\k $ be the class of $\a$ in $\I_\m/i(K_{\m,\Gamma})$, and suppose $\pi_0$ is a representation of $C^*(\a_\k\rtimes R_{\m,\Gamma}^*)$ on $H_0$.
If $\{\xi_n\}_n$ is an orthonormal basis for $\H_0$, then $\{u_{(x,1)}u_{(0,t_\a)}^*\otimes\xi_n : x\in\R_\a, n\geq 1\}$ is an orthonormal basis for $\H_\a$, so there is a natural unitary $\W_\a: \H_\a \to \ell^2(\R_\a) \otimes \H_0$
determined by $\W_{\tcb{\a}}: u_{(x,1)}u_{(0,t_\a)}^*\otimes\xi_n \mapsto \delta_{(x,\a)} \otimes \xi_n$ where $\{\delta_{(x,\a)}: x\in \R_\a\}$ is the standard basis for the finite-dimensional space $\ell^2(\R_{\tcb{\a}})$.
\end{lemma}
\begin{proof}
 For $(x,g)\in R\rtimes R_{\m,\Gamma}^*$ and $\xi\in\H_0$, we have
\[
u_{(x,g)}u_{(0,t_\a)}^*\otimes\xi=u_{(\bar{x}^\a,1)}u_{((x-\bar{x}^\a),g)}u_{(0,t_\a)}^*\otimes\xi=u_{(\bar{x}^\a,1)}u_{(0,t_\a)}^*\otimes \pi_0(u_{(t_\a (x-\bar{x}^\a),g)})\xi
\]
where the last equality uses that $t_\a (x-\bar{x}^\a)\in \a_\k$. This calculation implies that if $\{\xi_n\}$ is an orthonormal basis for $\H_0$, then the vectors $u_{(x,1)}u_{(0,t_\a)}^*\otimes\xi_n$ for $n\geq 1$ and $x\in\R_\a$ span a dense subspace of the Hilbert space $\X_\a\otimes_{C^*(\a_\k\rtimes R_{\m,\Gamma}^*)} \H_0$. Moreover, if $x,y\in\R_\a$ and $m,n \geq 1$, then 
\begin{align*}
\langle u_{(x,1)}u_{(0,t_\a)}^*\otimes\xi_n,u_{(y,1)}u_{(0,t_\a)}^*\otimes\xi_m\rangle&=\langle\pi_0(\langle u_{(y,1)}u_{(0,t_\a)}^*,u_{(x,1)}u_{(0,t_\a)}^*\rangle)\xi_n,\xi_m\rangle\\
&=\begin{cases}
\langle \pi_0(u_{(t_\a(x-y),1)})\xi_n,\xi_m\rangle & \text{ if } x-y\in\a, \\
0 & \text{ if } x-y\notin\a
\end{cases}\\
&=\begin{cases}
\langle \xi_n,\xi_m\rangle & \text{ if } x=y, \\
0 & \text{ if } x\neq y.
\end{cases}\\
&=\begin{cases}
 1 & \text{ if } x=y \text{ and } m=n \\
0 & \text{ if } x\neq y  \text{ or }m \neq n.
\end{cases}
\end{align*}
Hence  $\{u_{(x,1)}u_{(0,t_\a)}^*\otimes\xi_n : x\in\R_\a, n\geq 1\}$ is an orthonormal basis for $\H_\a$, and $\W_\a$ is the obvious unitary resulting from the natural bijection between the corresponding orthonormal bases.
\end{proof}

\begin{lemma}
\label{lem:Sb}
For each class $\k\in\I_\m/(K_{\m,\Gamma})$ define a Hilbert space by $\H_\k:=\bigoplus_{\a\in\k,\a\subseteq R} \H_\a$.
 Then there is an isometry $S_b :\H_\k \to \H_\k$ for each $b\in R_{\m,\Gamma}$ that is determined on the direct summand $\H_\a$ by
\begin{equation}
\label{eqn:isom}
S_b(zu_{(0,t_\a)}^*\otimes\xi)=u_{(0,b)}zu_{(0,b)}^*u_{(0,t_{b\a})}^*\otimes \pi_0(u_{(0,b\sofa^{-1}t_{b\a})})\xi
\qquad\quad
z\in C^*(R\rtimes R_{\m,\Gamma}^*), \  \xi\in\H_0.
\end{equation}
Moreover, $S_b$ maps $\H_{\tcb{\a}}$ into $\H_{\b\a}$.
\end{lemma}
\begin{proof} Let $b\in R_{\m,\Gamma}$.
Left multiplication by the element $u_{(0,b)}\in C^*(Q\rtimes K_{\m,\Gamma})$ defines a unitary
\[
\tilde{S}_b : C^*(Q\rtimes K_{\m,\Gamma})\otimes_{C^*(\a_\k\rtimes R_{\m,\Gamma}^*)}\H_0\to C^*(Q\rtimes K_{\m,\Gamma})\otimes_{C^*(\a_\k\rtimes R_{\m,\Gamma}^*)}\H_0
\]
such that $z\otimes\xi\mapsto u_{(0,b)}z\otimes\xi$ for all $z\in C^*(Q\rtimes K_{\m,\Gamma})$ and $\xi\in\H_0$.

For each integral ideal $\a\in\k$ and $b\in R_{\m,\Gamma}$, we have $(b^{-1}\sofa)b\a=\sofa\a=\a_\k$, which implies that $b\sofa^{-1}t_{b\a}\in R_{\m,\Gamma}^*$. Now if $z\in C^*(R\rtimes R_{\m,\Gamma}^*)$ and $\xi\in\H_0$, then
\[
u_{(0,b)}zu_{(0,\sofa)}^*\otimes\xi=u_{(0,b)}zu_{(0,b)}^*u_{(0,b^{-1}\sofa)}^*\otimes\xi=u_{(0,b)}zu_{(0,b)}^*u_{(0,t_{b\a})}^*\otimes \pi_0(u_{(0,b\sofa^{-1}t_{b\a})})\xi.
\]
Hence, restricting the unitary $\tilde{S}_b$ to $\H_\a\subseteq C^*(Q\rtimes K_{\m,\Gamma})\otimes_{C^*(\a_\k\rtimes R_{\m,\Gamma}^*)}\H_0$, we obtain an isometry $S_b^\a:\H_\a\to \H_{b\a}$ satisfying \eqref{eqn:isom}. Now $S_{\tcb{b}}:=\bigoplus_{\a\in\k, \a\subseteq R}S_b^\a$ has the desired properties.
\end{proof}

\begin{proposition}
\label{prop:extendrep} Fix a class $\k \in \I/i(K_{\m,\Gamma})$. For each 
 $x\in R$, let $U^x:=\bigoplus_{\a\in\k, \a\subseteq R}\vartheta_\a(u^x)$, and for each $b\in R_{\m,\Gamma}$, let $S_b$ be the isometry defined in \lemref{lem:Sb}. For integral ideals $\a,\b\in\I_\m$, let $\R_{\a/\b}:=\a\cap\R_\b$, so that $\R_{\a/\b}$ is a complete set of representatives for the (finite) quotient $\a/\b$, and $0\in\R_{\a/\b}$. 
Using the representative $\bar{x}^\a$ of $x$ in $\R_\a$, for $\a \in \k$ and $x \in R$, let
\[
\H_\a^x:=\overline{\spn}\{u_{(\bar{x}^\a,1)}u_{(0,\sofa)}^*\otimes\xi: \xi \in \H_0\}\subseteq \H_\a,
\] 
and define $E_\a$ to be the orthogonal projection from $\H_\k$ onto $\bigoplus_{\b\in\k, \a\mid\b, x\in\R_{\a/\b}} \H_\b^x$. Then there is a unique representation $\vartheta_\k$ of $C_\lambda^*(R\rtimes R_{\m,\Gamma})$ on 
$\H_\k:=\bigoplus_{\a\in\k,\a\subseteq R} \H_\a$ such that  $\vartheta_\k(u^xs_b)=U^xS_b$ for $(x,b)\in R\rtimes R_{\m,\Gamma}$. Moreover, $\vartheta_\k(e_\a)=E_\a$ for each integral ideal $\a\in\k$.
\end{proposition}
\begin{proof}
By \cite[Proposition~4.3]{Bru1}, it suffices to show that the families $\{U^x : x\in R\}$, $\{S_b: b\in R_{\m,\Gamma}\}$, and $\{E_\a : \a\in\I_\m , \a\subseteq R\}$ satisfy the following relations:
\begin{enumerate}
\item[(Ta)]  $U^xU^y=U^{x+y}$, $S_bS_c=S_{bc}$, and $S_bU^x=U^{bx}S_b$ for all $x,y\in R$ and $b,c\in R_{\m,\Gamma}$;

\smallskip
\item[(Tb)] $E_\a E_\b=E_{\a\cap \b}$  and $E_R=1$ for all integral $\a,\b\in\I_\m$;

\smallskip
\item[(Tc)] $S_bE_\a S_b^*=E_{b\a}$ for all $b\in R_{\m,\Gamma}$ and integral $\a\in\I_\m$;

\smallskip
\item[(Td)] $U^xE_\a=E_\a U^x$ if $x\in \a$ \ and \ $E_\a U^x E_\a =0$ if $x\not\in \a$.
\end{enumerate}
Proof of (Ta): Let $x,y\in R$ and $b,c\in R_{\m,\Gamma}$. We have $U^xU^y=U^{x+y}$ because each $\vartheta_\a$ is a *-homomorphism. 
Now $S_b$ and $S_c$ are defined via restriction to the $\H_\a$'s of left multiplication by $u_{(0,b)}$ and $u_{(0,c)}$, respectively, on $C^*(Q\rtimes K_{\m,\Gamma})\otimes_{C^*(\a_\k\rtimes R_{\m,\Gamma}^*)} \H_0$, and $U^x$ is defined via restriction to the $\H_\a$'s of left multiplication by $u_{(x,1)}$. Since $u_{(0,b)}u_{(0,c)}=u_{(0,bc)}$, we have $S_bS_c=S_{bc}$. Since $u_{(0,b)}u_{(x,1)}=u_{(bx,1)}u_{(0,b)}$, we have $S_bU^x=U^{bx}S_b$.

Proof of (Tb): For integral $\a,\b,\c\in\I_\m$, we have $\{\d\in\k :\d\subseteq R, \a\cap\b\mid \d \}=\{\d\in\k :\d\subseteq R, \a\mid \d \text{ and } \b\mid \d \}$ and $\R_{\a/\c}\cap\R_{\b/\c}=\R_{\a\cap\b/\c}$ whenever $\a\cap\b\mid \c$; it follows that $E_\a E_\b=E_{\a\cap \b}$. That $E_R=1$ follows directly from the definition of $E_R$.

Proof of (Tc): Let $b\in R_{\m,\Gamma}$ and suppose $\a\in\I_\m$ with $\a\subseteq R$. Let $\c$ be an integral ideal of $R$, $x\in\R_\c$, and $\xi\in\H_0$. Since the range of $S_b$ is equal to $\bigoplus_{\a\in\k, \a\subseteq R, y\in\R_\a}\H_{b\a}^{by}$, we have
\[
S_b^*(u_{(x,1)}u_{(0,t_\c)}^*\otimes\xi)=\begin{cases}
u_{(b^{-1}x,b^{-1})}u_{(0,t_\c)}^*\otimes\xi & \text{ if } x\in bR\text{ and } bR\mid \c ,\\
0 & \text{ otherwise.}
\end{cases}
\]
Moreover, for $x\in bR$ and $bR\mid \c$, we can write $S_b^*(u_{(x,1)}u_{(0,t_\c)}^*\otimes\xi)=u_{(\overline{b^{-1}x}^{b^{-1}\c},0)}u_{(0,t_{b^{-1}\c})}^*\otimes\eta$ for some $\eta\in\H_0$.
Hence, $E_\a S_b^*(u_{(x,1)}u_{(0,t_\c)}^*\otimes\xi)=0$ unless $x\in b\a$ and $b\a\mid \c$.

For $x\in\R_\c$, we have that $E_{b\a}(u_{(x,1)}u_{(0,t_\c)}^*\otimes\xi)$ is non-zero if and only if $x\in b\a$ and $b\a\mid \c$, and $S_bE_\a S_b^*$ is the identity on the range of $E_{b\a}$, so we have $S_bE_\a S_b^*=E_{b\a}$, as desired.

Proof of (Td): Let $\a\in\I_\m$ with $\a\subseteq R$ and $x\in R$. First, suppose $x\in\a$. To show $U^xE_\a=E_\a U^x$, it suffices to show that these operators agree on $\H_\b$ for $\b\in\I_\m$ with $\b\subseteq R$. Let $y\in\R_\b$. Since $x\in\a$, we can write $U^x(u_{(y,1)}u_{(0,t_\b)}^*\otimes\xi)=u_{(y,1)}u_{(0,t_\b)}^*\otimes\eta$ for some $\eta\in\H_0$ if $\a \mid \b$, so
\[
E_\a U^x(u_{(y,1)}u_{(0,t_\b)}^*\otimes\xi)
=\begin{cases}
u_{(y,1)}u_{(0,t_\b)}^*\otimes\eta & \text{ if } y\in\a \text{ and }\a\mid \b,\\
0 & \text{otherwise.}
\end{cases}
\]
On the other hand,
\[
U^x E_\a(u_{(y,1)}u_{(0,t_\b)}^*\otimes\xi)=\begin{cases}
U^x u_{(y,1)}u_{(0,t_\b)}^*\otimes\xi  = u_{(y,1)}u_{(0,t_\b)}^*\otimes\eta & \text{ if } y\in\a \text{ and }\a\mid \b,\\
0 & \text{otherwise}.
\end{cases}
\]
Hence $U^xE_\a=E_\a U^x$, as desired.

Now suppose $x\notin\a$. The range of $E_\a$ is $\bigoplus_{\a\in\k,\a\mid\b,y\in\R_{\a/\b}}\H_\b^y$. Since $x\notin\a$, $U^x\H_\b^y=\H_\b^{x+y}$ is orthogonal to $\H_\b^y$ for all $\b$. Hence, $E_\a U^xE_\a=0$, as desired.
\end{proof}

\begin{remark}\label{rem:Ea=UEaU}
The projection $E_\a^x:=U^xE_\a U^{-x}$ depends only on the class of $x$ modulo $\a$.
\end{remark}

In order to analyze this representation, it will be convenient to write $\bigoplus_\a \H_\a $ as 
$\ell^2\big(\bigsqcup _\a \R_\a\big) \otimes \H_0$.
\begin{lemma}
\label{lem:identifyrep}
Let $\k$ be a fixed class and suppose $\a$ is an integral ideal in $\k$.
Under conjugation by the unitary transformation 
\[
\W :  \ell^2\Big(\bigsqcup_{\a\in\k, \a\subseteq R}\R_\a\Big)\otimes \H_0 \overset{\cong}\longrightarrow \bigoplus_{\a\in\k, \a\subseteq R}\H_\a
\]
obtained by combining the unitaries $\W_\a$ from \lemref{lem:inducedrep},
the operators $U^x$, $S_b$, and $E_\a$ are given by
\begin{enumerate}
\item[(i)] $\W^*U^x \W(\delta_{(\c,y)}\otimes \xi)=\delta_{(\c,\overline{x+y}^\c)}\otimes\pi_0(u_{(t_\c(x+y-\overline{x+y}^\c),1)})\xi$;
\item[(ii)] $\W^*S_b \W(\delta_{(\c,y)}\otimes \xi)=\delta_{(b\c,\overline{by}^{b\c})}\otimes \pi_0(u_{(t_{b\c}(by-\overline{by}^{b\c}),bt_\c^{-1}t_{b\c})})\xi$;
\item[(iii)] the range of $\W^*E_\a \W$ is $\ell^2\big(\bigsqcup_{\b\in\k, \a\mid\b}\R_{\a/\b}\big)\otimes\H_0$;
\end{enumerate}
where $\a$ and $\c$ \tcb{are} integral ideals in the class $\k$, \ $x,y\in R$, \ $b\in R_{\m,\Gamma}$,  and $\xi\in\H_0$.
(After establishing this unitary equivalence, we will drop the unitary $\W$ from the formulas and abuse the  notation writing $U$, $S$, and $E$, for the families of operators  on $\ell^2\big(\bigsqcup_{\a\in\k, \a\subseteq R}\R_{\a}\big)\otimes\H_0$ defined by the right hand sides in (i)--(iii).)
\end{lemma}

\begin{proof}
Clearly, the $\W_\a$ from \lemref{lem:inducedrep} combine over the disjoint union, and their combined range is the direct sum, so $\W$ is a unitary transformation.
The explicit formula for $\W^* U^x \W$ in (i) follows from a short calculation, and the explicit formula  for $S_b$ in (ii) follows from \lemref{lem:Sb}. That the range of $\W^*E_\a\W$ is $\ell^2\big(\bigsqcup_{\b\in\k, \a\mid\b}\R_{\a/\b}\big)\otimes\H_0$ is immediate from the definition of $E_\a$.
\end{proof}

We are now ready to analyze the type of the representation of $C^*(\a_\k\rtimes R_{\m,\Gamma}^*)$ generated by $U$, $S$ and $E$ in terms of the type of the representation $\pi_0$. \tcb{We will rely on \cite{Bla} as a general reference for the necessary facts about von Neumann algebras.}

\begin{theorem}
\label{thm:vNalg}
Fix a class $\k$ in $\I_\m/i(K_{\m,\Gamma})$, suppose $\pi_0$ is a representation of $C^*(\a_\k\rtimes R_{\m,\Gamma}^*)$ on $\H_0$, and let $\vartheta_\k$ be the associated representation of $C_\lambda^*(R\rtimes R_{\m,\Gamma})$ on $\ell^2\big(\bigsqcup_{\a\in\k, \a\subseteq R}\R_\a\big)\otimes \H_0 = \bigoplus_{\a\in\k, \a\subseteq R} \ell^2(\R_\a; H_0)$ from \proref{prop:extendrep}~ via \lemref{lem:identifyrep}. 
Denote the associated von-Neumann algebras by  $\M:=\vartheta_\k(C_\lambda^*(R\rtimes R_{\m,\Gamma}))''$ and $\N:=\pi_0(C^*(\a_\k\rtimes R_{\m,\Gamma}^*))''\subseteq \B(\H_0)$. Then
\begin{equation}\label{eqn:vNtensorproduct}
\M=\B\Big(\ell^2\Big(\bigsqcup_{\a\in\k, \a\subseteq R}\R_\a\Big)\Big){\bar\otimes}\N
\end{equation}
under the natural identification of $ \B\big(\ell^2\big(\bigsqcup_{\a\in\k, \a\subseteq R}\R_\a\big) \otimes\H_0 \big)= \B\big(\ell^2\big(\bigsqcup_{\a\in\k, \a\subseteq R}\R_\a\big)\big){\bar\otimes}\B(\H_0)$.
Hence,
\begin{itemize}
\item[(i)] if $\N$ is a type I (resp. type II$_1$) factor, then $\M$ is a type I$_\infty$ (resp. type II$_\infty$) factor;
\item[(ii)] $\vartheta_\k$ is irreducible if and only if $\pi_0$ is irreducible. 
\end{itemize} 
\end{theorem}

\begin{proof}[Proof of Theorem~\ref{thm:vNalg}]
To prove \eqref{eqn:vNtensorproduct} it suffices to show that we have the following three inclusions:
\begin{enumerate}
\item[(a)] $\M\supseteq \B\big(\ell^2\big(\bigsqcup_{\a\in\k, \a\subseteq R}\R_\a\big)\big)\bar{\otimes}\CC$;
\item[(b)] $\M\supseteq \CC\bar{\otimes}\N$;
\item[(c)] $\M\subseteq \B\big(\ell^2\big(\bigsqcup_{\a\in\k, \a\subseteq R}\R_\a\big)\big)\bar{\otimes}\N$. 
\end{enumerate}

We introduce first some notation that will be useful in several places. 
Let $e_{R\times \a^\times}:=\sum_{x\in R/\a} e_{x+\a}$. Then the operator $\vartheta_\k(e_{R\times \a^\times})=\sum_{x\in\R_\a}\vartheta_\k(u^xe_\a u^{-x})=\sum_{x\in\R_\a}U^xE_\a U^{-x}$ is the orthogonal projection onto the subspace 
\[
\ell^2\Big(\bigsqcup_{\b\in\k, \a\mid \b} \R_\b\Big)\otimes\H_0.
\]
Therefore, the projection $\E_\a:=\vartheta_\k(e_{R\times \a^\times})-\bigvee_{\b\in\k, \b\subsetneq \a}\vartheta_\k(e_{R\times \b^\times})$ lies in $\M$. Note that $\E_\a$ is the orthogonal projection onto the subspace $\ell^2(\R_\a)\otimes\H_0\cong \H_\a$, and, recalling the projection $E_\a^x:=U^xE_\a U^{-x}$ from \remref{rem:Ea=UEaU},
observe that $E_\a^x$ and $\E_\a$ commute and that $\E_\a E_\a^x$ is the orthogonal projection onto $\CC\delta_{(\a,x)}\otimes \H_0$. 

Proof of (a): For integral ideals $\a,\b\in\k$, $x\in \R_\a$, and $y\in \R_\b$, let $W_{\a,\b}^{x,y}$ be the partial isometry in $\B\big(\ell^2\big(\bigsqcup_{\a\in\k, \a\subseteq R}\R_\a\big)\big)$ with initial space $\CC\delta_{(\a,x)}$ and final space $\CC\delta_{(\b,y)}$, so that $(W_{\a,\b}^{x,y})^*=W_{\b,\a}^{y,x}$ and $(W_{\a,\b}^{x,y}\otimes 1)^*(W_{\a,\b}^{x,y}\otimes 1)=\E_\a E_\a^x$.
It is enough to show that each $W_{\a,\b}^{x,y}\otimes 1$ belongs to $\M$. First, suppose that there exists $b\in R_{\m,\Gamma}$ with $\b=b\a$. 
Since $t_{b\a}^{-1}t_\a b^{-1}\in R_{\m,\Gamma}^*$, we have $t_{b\a}^{-1}t_\a\in R_{\m,\Gamma}$. Moreover, note that $t_{b\a}^{-1}t_\a\a=t_{b\a}^{-1}\a_\k=b\a=\b$. The operator $U^yS_{t_{b\a}^{-1}t_\a}U^{x*}\E_\a E_\a^x$ is zero on $(\CC\delta_{(\a,x)}\otimes \H_0)^\perp\subseteq \ell^2\big(\bigsqcup_{\a\in\k, \a\subseteq R}\R_\a\big)\otimes\H_0$, and for every $\xi\in\H_0$, we can compute using Lemma~\ref{lem:identifyrep}:
\begin{align*}
U^yS_{t_{b\a}^{-1}t_\a}U^{x*}(\delta_{(\a,x)} \otimes \xi)=U^yS_{t_{b\a}^{-1}t_\a}(\delta_{(\a,0)} \otimes\xi)&=U^y(\delta_{(t_{b\a}^{-1}t_\a\a,0)} \otimes \pi_0(u_{(0,(t_\b^{-1}t_\a)t_\a^{-1}t_\b)})\xi)\\
&=\delta_{(\b,y)} \otimes\xi 
\end{align*}
where we used that $\bar{x}^\a=x$ and $\bar{y}^\b=y$ because \tcb{$x\in \R_\a$} and $y\in \R_\b$. Thus, $W_{\a,\b}^{x,y}\otimes 1=U^yS_{t_{b\a}^{-1}t_\a}U^{x*}\E_\a E_\a^x\in\M$.

For the general case, we define an equivalence relation on the integral ideals in $\k$ as follows: for $\a,\b\in\k$ with $\a,\b\subseteq R$, we say $\a\sim \b$ if and only if there exists $x\in\R_\a$ and $y\in\R_\b$ such that $W_{\a,\b}^{x,y}\otimes 1\in\M$. The above argument for the case $b=1$ shows that $\sim$ is reflexive. Symmetry follows from $(W_{\a,\b}^{x,y}\otimes 1)^*=W_{\b,\a}^{y,x}\otimes 1$, and transitivity follows from taking products. Hence,  $\sim$ is an equivalence relation.
Since $\a,\b\in\k$, there are $c,d\in R_{\m,\Gamma}$ such that $c\a=d\b$. Now
\[
\a\sim c\a\sim d\b\sim \b,
\] 
as desired.

Proof of (b): It suffices to show that $\M$ contains the operators $W_{\a,\a}^{y,y}\otimes \pi_0(u_{(x,1)})$ and $W_{\a,\a}^{y,y}\otimes \pi_0(u_{(0,g)})$ for each fixed integral ideal $\a\in\k$, $x\in \a_\k$, $y\in\R_\a$, and $g\in R_{\m,\Gamma}^*$ where $W_{\a,\a}^{y,y}$ is the projection in $\B\big(\ell^2\big(\bigsqcup_{\a\in\k, \a\subseteq R}\R_\a\big)\big)$ onto the subspace $\CC\delta_{(\a,y)}$.

Since $\a=\sofa^{-1}\a_\k$, we have $\sofa^{-1}x \in \a$. For every $\xi\in\H_0$, Lemma~\ref{lem:identifyrep} gives us
\[
U^{\sofa^{-1}x}(\delta_{(\a,y)}\otimes\xi)=\delta_{(\a,\overline{\sofa^{-1}x+y}^\a)}\otimes \pi_0(u_{(t_\a(\sofa^{-1}x+y-\overline{\sofa^{-1}x+y}^\a),1)})\xi=\delta_{(\a,y)}\otimes\pi_0(u_{(x,1)})\xi.
\]
where we used that $\overline{\sofa^{-1}x+y}^\a=y$ (which holds because $\sofa^{-1}x\in\a$ and $y\in\R_\a$).
Hence, $W_{\a,\a}^{y,y}\otimes \pi_0(u_{(x,1)})=U^{\sofa^{-1}x}\E_\a E_\a^y$, which lies in $\M$. 

Next, let $c:=\overline{gy}-gy$, where $\overline{gy}$ is the unique elements in $\R_\a$ such that $\overline{gy}-gy\in \a$. Let $\xi\in \H_0$. Since $\overline{c+\overline{gy}}=\overline{gy}$, Lemma~\ref{lem:identifyrep} yields
\begin{align*}
U^cS_g(\delta_{(\a,y)}\otimes\xi)&=U^c(\delta_{(\a,\overline{gy})}\otimes\pi_0(u_{(t_\a(gy-\overline{gy}),g)})\xi=\delta_{(\a,\overline{gy})}\otimes \pi_0(u_{(t_\a(c+\overline{gy}-\overline{c+\overline{gy}}),1)})\pi_0(u_{(t_\a(gy-\overline{gy}),g)})\xi\\
&=\delta_{(\a,\overline{gy})}\otimes \pi_0(u_{(0,g)})\xi,
\end{align*}
so $W_{\a,\a}^{y,\overline{gy}}\otimes \pi_0(u_{(0,g)})=U^cS_g\E_\a E_\a^y$, which lies in $\M$. Hence,
\[
W_{\a,\a}^{y,y}\otimes \pi_0(u_{(0,g)})=(W_{\a,\a}^{y,\overline{gy}}\otimes 1)^*(W_{\a,\a}^{y,\overline{gy}}\otimes\pi_0(u_{(0,g)}))
\]
also lies in $\M$.

Proof of (c): Since elements of the form $U^xS_b=\vartheta_\k(u^xs_b)$ for $(x,b)\in R\rtimes R_{\m,\Gamma}$ generate a dense *-subalgebra of $\vartheta_\k(C^*_\lambda(R\rtimes R_{\m,\Gamma}))$, it suffices to show that $U^x$ and $S_b$ are in $\B\big(\ell^2\big(\bigsqcup_{\a\in\k, \a\subseteq R}\R_\a\big)\big)\bar{\otimes}\N$ for all $x\in R$ and $b\in R_{\m,\Gamma}$. This follows from the explicit formulas for $U^x$ and $S_b$ given in Lemma~\ref{lem:identifyrep}. Indeed, given an integral ideal $\a\in\I_\m$, $y\in\R_\a$, and $\xi\in\H_0$, we have $U^x(\delta_{(\a,y)}\otimes\xi)=\delta_{(\a,\overline{x+y}^\a)}\otimes\pi_0(z)\xi$ and $S_b(\delta_{(\a,y)}\otimes\xi)=\delta_{(b\a,\overline{by}^{b\a})}\otimes\pi_0(z')\xi$ for some $z,z'\in C^*(\a_\k\rtimes R_{\m,\Gamma}^*)$. Hence the operators $U^x$ and $S_b$ are given by infinite matrices with entries in $\N$, and are thus in $\B\big(\ell^2\big(\bigsqcup_{\a\in\k, \a\subseteq R}\R_\a\big)\big)\bar{\otimes}\N$.

Once we have the factorization \eqref{eqn:vNtensorproduct} it is easy to see that  assertions (i) and  (ii) hold,
the latter because of the standard fact $\left(\B\big(\ell^2\big(\bigsqcup_{\a\in\k, \a\subseteq R}\R_\a\big)\big)\bar{\otimes}\N\right)'=\CC I\bar{\otimes} \N'$, \tcb{see, for instance, \cite[III.1.5.9]{Bla}}.
\end{proof}

The following result, which explicitly constructs KMS states on $C_\lambda^*(R\rtimes R_{\m,\Gamma})$ from tracial states on $\N$, will allow us to use Theorem~\ref{thm:vNalg} to both determine type and compute partition functions.


\begin{proposition}
\label{prop:tracestostates}
Resume the notation and the assumptions from Theorem~\ref{thm:vNalg}, and for each class $\k$ 
let $\unk$ be the smallest among the norms of the integral ideals in $\k$.
For each $t\in\RR$, let $e^{itH_\k}$ denote the diagonal unitary operator on $\ell^2\big(\bigsqcup_{\a\in\k, \a\subseteq R}\R_\a\big)$ such that $e^{itH_\k}(\delta_{(\a,x)})=\left(\frac{N(\a)}{\unk}\right)^{it}\delta_{(\a,x)}$.
Then $t\mapsto \Ad(e^{itH_\k}\otimes I)$ is a time evolution on $\ell^2\big(\bigsqcup_{\a\in\k, \a\subseteq R}\R_\a\big)\otimes\H_0$, and we have 
\begin{equation}
\label{eqn:equivariant}
\vartheta_\k(\sigma_t(z))=\Ad(e^{itH_\k}\otimes I)(\vartheta_\k(z))\quad \text{for }t\in\RR,\;z\in C_\lambda^*(R\rtimes R_{\m,\Gamma}).
\end{equation}

For each $\beta\in (2,\infty)$, the diagonal operator $e^{-\beta H_\k}$ defined by $e^{-\beta H_\k}(\delta_{(\a,x)})=N(\a)^{-\beta}(\delta_{(\a,x)})$ is of trace class, and $\Tr(e^{-\beta H_\k})=\unk^{\tcb{\beta}}\zeta_\k(\beta-1)$; let $\psi_{\beta,\k}$ denote the Gibbs state on $\M$ defined by $\psi_{\beta,\k}(T):=\frac{\Tr(Te^{-\beta H_\k})}{\unk^{\tcb{\beta}}\zeta_\k(\beta-1)}$.
If $\tilde{\tau}$ is a normal tracial state on $\N$, then 
\[
\psi_{\beta,\k,\tilde{\tau}}:=(\psi_{\beta,\k}\otimes\tilde{\tau})\circ \vartheta_\k=(\Tr\otimes\tilde{\tau})\Big(\vartheta_\k(\cdot) \frac{e^{-\beta H_\k}}{\unk^{\tcb{\beta}}\zeta_\k(\beta-1)}\otimes I\Big)
\] is a $\sigma$-KMS$_\beta$ state on $C_\lambda^*(R\rtimes R_{\m,\Gamma})$,
 and  if we let   $\tau:=\tilde{\tau}\circ\pi_0$, then  $\psi_{\beta,\k,\tilde{\tau}}$ is equal to the state $\varphi_{\beta,\k,\tau}$ from \proref{pro:KMSformula}.
\end{proposition}
\begin{proof}
It is routine to check that $t\mapsto \Ad(e^{itH_\k}\otimes I)$ is a time evolution on $\ell^2\big(\bigsqcup_{\a\in\k, \a\subseteq R}\R_\a\big)\otimes\H_0$. To show that $\vartheta_\k$ is $\RR$-equivariant, it suffices to show that \eqref{eqn:equivariant} holds for $z$ of the form $u^x$ or $s_b$ where $x\in R$, $b\in R_{\m,\Gamma}$. This is easy to see for $z=u^x$. Take an integral ideal $\a\in\k$, $x\in\R_\a$, and $\xi\in\H_0$. For $t\in\RR$,  Lemma~\ref{lem:identifyrep}(ii) implies that there exist $\eta\in\H_0$ such that $S_b=\delta_{(b\a,\overline{bx}^{b\a})}\otimes\eta$. Hence,
\begin{align*}
\Ad(e^{itH_\k}\otimes I)(\vartheta_\k(s_b))\delta_{(\a,x)}\otimes\xi &=(e^{itH_\k}\otimes I)S_b(e^{-itH_\k}\otimes I)\delta_{(\a,x)}\otimes\xi=\left(\frac{N(\a)}{\unk}\right)^{-it}(e^{itH_\k}\otimes I)\delta_{(b\a,\overline{bx}^{b\a})}\otimes\eta\\
&=N(b)^{it}\delta_{(b\a,\overline{bx}^{b\a})}\otimes\eta=N(b)^{it}S_b\delta_{(\a,x)}\otimes\xi=\vartheta_\k(\sigma_t(s_b))\delta_{(\a,x)}\otimes\xi,
\end{align*}
from which $\RR$-equivariance follows.

A short calculation using the orthonormal basis $\{\delta_{(\a,x)} : \a\in\I_\m,\a\subseteq R, x\in\R_\a\}$ for $\ell^2\big(\bigsqcup_{\a\in\k, \a\subseteq R}\R_\a\big)$ shows that $e^{-\beta H_\k}$ is of trace class with $\Tr(e^{-\beta H_\k})=\unk^{\beta}\zeta_\k(\beta-1)$ if $\beta > 2$.

Lastly, we must show that $\psi_{\beta,\k,\tilde{\tau}}=\varphi_{\beta,\k,\tau}$ where $\tau:=\tilde{\tau}\circ\pi_0$. For this, we first prove that it is enough to deal with the case where $\tilde{\tau}$ is a vector state. 
Assume the result holds in this case, and let $\tilde{\tau}$ any tracial state on $\N$. Let $\pi_0':C^*(\a_\k\rtimes R_{\m,\Gamma}^*)\to \B(\H_0\otimes\H_0)$ by $\pi_0'(z)=\pi_0(z)\otimes I$. Since $\tilde{\tau}$ is normal, there exists orthogonal vectors $\xi_1,\xi_2,...$ in $\H_0$ such that $\tilde{\tau}(\cdot)=\sum_n\langle\cdot\xi_n,\xi_n\rangle$ and $\sum_n \|\xi_n\|^2=1$ (see, for instance, \cite[Theorem~III.2.1.4(vii)]{Bla}). Then a calculation shows that $\eta:=\sum_n\xi_n\otimes\xi_n\in\H_0\otimes\H_0$ is a unit vector and the corresponding vector state $\tilde{\tau}'$ satisfies $\tilde{\tau}'(z\otimes I)=\tilde{\tau}(z)$ for all $z\in \N$, so that $\tilde{\tau}'\circ\pi_0'=\tilde{\tau}\circ\pi_0=\tau$.
Let $\vartheta_\k'$ denote the representation associated with $\pi_0'$ using Proposition~\ref{prop:extendrep}. By assumption, $\psi_{\beta,\k,\tilde{\tau}'}=\varphi_{\beta,\k,\tau}$. Thus, it remains to show that $\psi_{\beta,\k,\tilde{\tau}}=\psi_{\beta,\k,\tilde{\tau}'}$. For this, it suffices to show that the following diagram commutes:
\[\begin{tikzcd}
 & \M\otimes I\subseteq\B\big(\ell^2\big(\bigsqcup_{\a\in\k, \a\subseteq R}\R_\a\big)\otimes\H_0\otimes\H_0\big)\\
C_\lambda^*(R\rtimes R_{\m,\Gamma}) \arrow[r,"\vartheta_\k"]\arrow[ru,"\vartheta_\k'"] & \M\subseteq \B\big(\ell^2\big(\bigsqcup_{\a\in\k, \a\subseteq R}\R_\a\big)\otimes\H_0\big).\arrow[u,"T\mapsto T\otimes I"']
\end{tikzcd}\]
Commutativity follows from the definitions of $\vartheta_\k$ and $\vartheta_\k'$ combined with the fact that the canonical isomorphism of Hilbert spaces
\[
C^*(Q\rtimes K_{\m,\Gamma})\otimes_{\pi_0'}(\H_0\otimes\H_0)\cong (C^*(Q\rtimes K_{\m,\Gamma})\otimes_{\pi_0}\H_0)\otimes\H_0
\]
is even an isomorphism of left $C^*(Q\rtimes K_{\m,\Gamma})$-modules. This concludes our proof that it suffices to deal with the case where $\tilde{\tau}$ is a vector state.

Now suppose that $\tilde{\tau}(\cdot)=\langle\cdot\eta,\eta\rangle$ for some unit vector $\eta\in \H_0$. Then a short calculation shows that
\[
(\psi_{\beta,\k}\otimes\tilde{\tau})(\cdot)=\frac{1}{\zeta_\k(\beta-1)}\sum_{\a\in\k,\a\subseteq R,x\in\R_\a}N(\a)^{-\beta}\langle\cdot \delta_{(\a,x)}\otimes\eta,\delta_{(\a,x)}\otimes\eta\rangle.
\]

Therefore, we we have
\[
\psi_{\beta,\k,\tilde{\tau}}(s_b^*e_{y+\b}u^ds_c)=(\psi_{\beta,\k}\otimes\tilde{\tau})(S_b^*E_\b^yU^dS_c)=\frac{1}{\zeta_\k(\beta-1)}\sum_{\a\in\k,\a\subseteq R,x\in\R_\a}N(\a)^{-\beta}\langle S_b^*E_\b^yU^dS_c \delta_{(\a,x)}\otimes\eta,\delta_{(\a,x)}\otimes\eta\rangle.
\]
Using Proposition~\ref{lem:identifyrep}(i)\&(ii), we have
\begin{enumerate}
\item $S_b(\delta_{(\a,x)}\otimes\eta)=\delta_{(b\a,\overline{bx}^{b\a})}\otimes\pi_0(u_{(t_{b\a}(bx-\overline{bx}^{b\a}),bt_\a^{-1}t_{b\a})})\eta$;
\item $U^dS_c(\delta_{(\a,x)}\otimes\eta)=\delta_{(c\a,\overline{d+\overline{cx}}^{c\a})}\otimes\pi_0(u_{(t_{c\a}(d+cx-\overline{d+\overline{cx}}^{c\a}),ct_\a^{-1}t_{c\a})})\eta$, where $\overline{cx}$ is the unique element in $\R_{c\a}$ such that $cx-\overline{cx}\in c\a$
\end{enumerate} 

Now, for
\[\langle S_b^*E_\b^yU^dS_c \delta_{(\a,x)}\otimes\eta,\delta_{(\a,x)}\otimes\eta\rangle=\langle E_\b^yU^dS_c \delta_{(\a,x)}\otimes\eta,S_b\delta_{(\a,x)}\otimes\eta\rangle
\]
to be non-zero, we must have $b\a=c\a$, $\overline{bx}^{c\a}=\overline{d+\overline{cx}}^{c\a}$, $\b\mid c\a$, and $\overline{d+\overline{cx}}^{c\a}=\overline{y+z}^{c\a}$ for some $z\in \R_{\b/c\a}$. The second condition is equivalent to $d+cx-bx\in c\a$, and the fourth condition is equivalent to $d+cx-y\in\b$. When these four conditions are satisfied, we have 
\begin{align*}
\langle S_b^*E_\b^yU^dS_c \delta_{(\a,x)}\otimes\eta,\delta_{(\a,x)}\otimes\eta\rangle&=\langle\pi_0(u_{(t_{c\a}(d+cx-\overline{d+\overline{cx}}^{c\a}),ct_\a^{-1}t_{c\a})})\eta,\pi_0(u_{(t_{b\a}(bx-\overline{bx}^{c\a}),bt_\a^{-1}t_{b\a})})\eta\rangle\\
&=\langle\pi_0(u_{(t_{b\a}(bx-\overline{bx}^{c\a}),bt_\a^{-1}t_{b\a})})^*\pi_0(u_{(t_{c\a}(d+cx-\overline{d+\overline{cx}}^{c\a}),ct_\a^{-1}t_{c\a})})\eta,\eta\rangle\\
&=\langle\pi_0(u_{(b^{-1}t_\a(d+cx-bx),b^{-1}c)})\eta,\eta\rangle=\tau(u_{(b^{-1}t_\a(d+cx-bx),b^{-1}c)}).
\end{align*}

Thus,
\begin{align*}
\frac{1}{\zeta_\k(\beta-1)}\sum_{\substack{\a\in\k,\a\subseteq R,\\x\in\R_\a}}N(\a)^{-\beta}&\langle S_b^*E_\b^yU^dS_c \delta_{(\a,x)}\otimes\eta,\delta_{(\a,x)}\otimes\eta\rangle\\
&= \begin{cases} \displaystyle \frac{1}{\zeta_\k(\beta-1)}\sum_{\substack{\a\in\k, \b\mid c\a,\\ \a \subseteq R}}\sum_{\substack{x\in R/\a,\\d+cx-y\in\b,\\ d+cx-bx\in c\a}}N(\a)^{-\beta}\tau(u_{(b^{-1}t_\a(d+cx-bx),b^{-1}c)}) & \text{ if } b^{-1}c \in R_{\m,\Gamma}^* \\
\quad 0 & \text{ otherwise} \end{cases}
\end{align*}
which is precisely the explicit formula for $\varphi_{\beta,\k,\tau}$ given in Proposition~\ref{pro:KMSformula}.
\end{proof}

\begin{remark}
If $\pi_0=\bigoplus_{\varphi\in S(C^*(\a_\k\rtimes R_{\m,\Gamma}^*))}\pi_\varphi$ is the universal representation of $C^*(\a_\k\rtimes R_{\m,\Gamma}^*)$, so that $\N=\pi_0(C^*(\a_\k\rtimes R_{\m,\Gamma}^*))''\cong C^*(\a_\k\rtimes R_{\m,\Gamma}^*)^{**}$ and each tracial state $\tau\in T(C^*(\a_\k\rtimes R_{\m,\Gamma}^*))$ has a unique extension to a tracial state $\tilde{\tau}$ on $\N$, then 
\[
T(C^*(\a_\k\rtimes R_{\m,\Gamma}^*))\ni\tau \mapsto (\psi_{\beta,\k}\otimes\tilde{\tau})\circ \vartheta_\k\in \textup{KMS}_\beta(C_\lambda^*(R\rtimes R_{\m,\Gamma}),\sigma)
\]
defines an embedding of simplices, and these embeddings combine to give an isomorphism of simplices
\[
T\big(\bigoplus_{\k\in\I_\m/i(K_{\m,\Gamma})}C^*(\a_\k\rtimes R_{\m,\Gamma}^*)\big)\cong \textup{KMS}_\beta(C_\lambda^*(R\rtimes R_{\m,\Gamma}),\sigma)
\]
(cf. \cite[Theorem~3.2(iii)]{Bru2}).

However, in our main applications, $\N$ will have a unique normal tracial state arising from a tracial state on $C^*(\a_\k\rtimes R_{\m,\Gamma}^*)$.
\end{remark}

We are now ready for our computation of type.
In our main application, the representation $\pi_0$  will be an irreducible tensor factor of the GNS representation  of an extremal type I tracial state of $C^*(\a_\k\rtimes R_{\m,\Gamma}^*)$ or the GNS representation of an extremal type II tracial state of $C^*(\a_\k\rtimes R_{\m,\Gamma}^*)$. 
\begin{corollary}
\label{cor:typeKMS}
Let $\k\in \I_\m/i(K_{\m,\Gamma})$,  fix $\beta\in (2,\infty)$, and let $\tau$ 
be an extremal trace on  $C^*(\a_\k\rtimes R_{\m,\Gamma}^*)$.
If $\tau$ is type I, then $\varphi_{\beta,\k,\tau}$ is type I$_\infty$; and 
if $\tau$ is type II$_1$, then $\varphi_{\beta,\k,\tau}$ is type II$_\infty$.
By Proposition~\ref{prop:typeoftraces}  this exhausts all possibilities.
\end{corollary}
\begin{proof}
Let $\pi_0=\pi_\tau$ be the GNS representation associated with the tracial state $\tau\in T(C^*(\a_\k\rtimes R_{\m,\Gamma}^*))$, so that $\N=\pi_0(C^*(\a_\k\rtimes R_{\m,\Gamma}^*))''\subseteq \B(\H_\tau)$. Since 
\tcb{$\varphi:=\varphi_{\beta,\k,\tau}$} is extremal, $\pi_\varphi(C_\lambda^*(R\rtimes R_{\m,\Gamma}))''$ is a factor. Since $\tau$ is extremal, so that $\N$ is a factor, Theorem~\ref{thm:vNalg} implies that $\M$ is also a factor of type I$_\infty$ or type II$_\infty$ according to whether $\N$ is of type I or type II. Thus, it is enough to show that $\pi_\varphi(C_\lambda^*(R\rtimes R_{\m,\Gamma}))''$ is isomorphic to $\M$.

By Proposition~\ref{prop:tracestostates}, the KMS state $\varphi \tcb{=}\varphi_{\beta,\k,\tau}$ is given explicitly by 
\[
\varphi(z)=(\Tr\otimes\tilde{\tau})\Big(\vartheta_\k(z) \frac{e^{-\beta H_\k}}{\unk^{\tcb{\beta}}\zeta_\k(\beta-1)}\otimes I\Big)
\] 
for all $z\in C_\lambda^*(R\rtimes R_{\m,\Gamma})$. Let $(\pi_\psi,\H_\psi,\xi_\psi)$ be the GNS representation of $\M$ associated with the normal state $\psi:x\mapsto (\Tr\otimes\tilde{\tau})\big(x\big(\frac{e^{-\beta H_\k}}{\unk^{\tcb{\beta}}\zeta_\k(\beta-1)}\otimes I\big)\big)$ of $\M$. We have $\varphi=\psi\circ\vartheta_\k$, so
\[
\langle\pi_\psi\circ\vartheta_\k(z)\xi_\psi,\xi_\psi\rangle=\psi(\vartheta_\k(z))=\varphi(z)
\]
for all $z\in C_\lambda^*(R\rtimes R_{\m,\Gamma})$. Now, by uniqueness of the GNS representation, $\pi_{\varphi}$ is unitarily equivalent to $\pi_\psi\circ\vartheta_\k$. Since $\pi_\varphi(C_\lambda^*(R\rtimes R_{\m,\Gamma}))''$ and $\M$ are factors, we have isomorphisms $\pi_\varphi(C_\lambda^*(R\rtimes R_{\m,\Gamma}))''\cong \pi_\psi(\M)\cong\M$.
\end{proof}

Next we compute  an admissible triple and the partition function for each extremal KMS$_\beta$ state on $C_\lambda^*(R\rtimes R_{\m,\Gamma})$ that is of type I.

\begin{theorem}
\label{thm:partitionfncI}
Let $\varphi=\varphi_{\k,\beta_0,\tau}$ be the extremal KMS$_{\beta_0}$ state on $C_\lambda^*(R\rtimes R_{\m,\Gamma})$ corresponding to a class $\k\in \I_\m / i(K_{\m,\Gamma})$, an inverse temperature $\beta_0\in(2,\infty)$, and an extremal trace $\tau=\tau_{O,\chi}$ on $C^*(\a_\k\rtimes R_{\m,\Gamma}^*)$ arising from the finite orbit $O$ and the character $\chi$ of its isotropy subgroup $(R_{\m,\Gamma}^*)_O$. Then the partition function of $\varphi$  from Definition~\ref{def:partitionfunction} is given by
\[
Z_\varphi(s)=|O|\, \unk^s \, \zeta_\k(s-1)\quad \text{for } \Re s > 2.
\]

Moreover, the irreducible representations of admissible triples and the GNS representation of $\varphi$ are faithful.
\end{theorem}
\begin{proof}
For each $x\in O$, let $\H_{x,\chi}:=C^*(R_{\m,\Gamma}^*)\otimes_{C^*((R_{\m,\Gamma}^*)_O)} \CC$ where $C^*((R_{\m,\Gamma}^*)_O)$ acts on $\CC$ via $\chi$. Then there is an $|O|$-dimensional irreducible representation $\pi_{x,\chi}$ of $C^*(\a_\k\rtimes R_{\m,\Gamma}^*)\cong C(\widehat{\a}_\k)\rtimes R_{\m,\Gamma}^*$ on $\H_{x,\chi}$ given by 
\[
\pi_{x,\chi}(f)(u_g\otimes 1)=f(gx)u_g\otimes 1 \text{ and } \pi_{x,\chi}(u_h)(u_g\otimes 1)=u_{gh}\otimes 1,
\]
where $f\in C(\widehat{\a}_\k)$ and $g,h\in R_{\m,\Gamma}^*$. Hence, $\N=\pi_{x,\chi}(C^*(\a_\k\rtimes R_{\m,\Gamma}^*))''=\pi_{x,\chi}(C^*(\a_\k\rtimes R_{\m,\Gamma}^*))=\B(\H_{x,\chi})$.
Let $\vartheta_\k$ be the representation of $C_\lambda^*(R\rtimes R_{\m,\Gamma})$ associated with $\pi_0:=\pi_{x,\chi}$ from Proposition~\ref{prop:extendrep}. By Theorem~\ref{thm:vNalg}(ii), $\vartheta_\k$ is irreducible. 

By Proposition~\ref{prop:tracestostates}, 
\[
\varphi=(\Tr_{\ell^2}\otimes\frac{1}{|O|}\Tr_{\H_{x,\chi}})\Big(\vartheta_\k(\cdot) \frac{e^{-\beta_0 H_\k}}{\unk^{\beta_0}\zeta_\k(\beta_0-1)}\otimes I\Big)=(\Tr_{\ell^2}\otimes\Tr_{\H_{x,\chi}})\Big(\frac{\vartheta_\k(\cdot)(e^{-\beta_0 H_\k}\otimes I)}{|O| \, \unk^{\beta_0} \, \zeta_\k(\beta_0-1)}\Big)
\] 
where $\Tr_{\ell^2}$ and $\Tr_{\H_{x,\chi}}$ are the canonical (unnormalized) traces on $\B\big(\ell^2\big(\bigsqcup_{\a\in\k, \a\subseteq R}\R_\a\big)\big)$ and $\B(\H_{x,\chi})$, respectively. Now, $\Tr_{\ell^2}\otimes\Tr_{\H_{x,\chi}}$ is the canonical trace on $\B\big(\ell^2\big(\bigsqcup_{\a\in\k, \a\subseteq R}\R_\a\big)\otimes\H_{x,\chi}\big)$, and the operator $e^{-\beta_0 H_\k}\otimes I$ has norm $1$ and is of trace-class with trace $|O| \, \unk^{\beta_0} \, \zeta_\k(\beta_0-1)$. Hence the triple consisting of the representation $\vartheta_\k$, the Hilbert space  $\ell^2\big(\bigsqcup_{\a\in\k, \a\subseteq R}\R_\a\big)\otimes\H_{x,\chi}$ and the density operator $e^{-\beta_0 H_\k}\otimes I$ is admissible for $\varphi$ so the partition function $Z_\varphi(\beta)$ is equal to 
$|O| \, \unk^s \, \zeta_\k(\beta-1)$ (see Definition~\ref{def:partitionfunction}).

By \proref{pro:gnsfactorization}, to prove the last assertion it suffices to show that the GNS representation $\pi_\varphi$ of $\varphi$ is
 faithful. For this, we shall use \cite[Theorem~6.1]{Bru1}. 
 Recall from the proof of  \proref{pro:KMSfromtrace} that the restriction  $\varphi\vert_{C(\Omega)}$ of $\varphi$ to the diagonal $C(\Omega)$
 is given by integration against the probability measure $\mu_{\beta,\k}$ given by \eqref{eqn:measures}.
Suppose $\tilde{\k}$ is a class in $ \I_\m/i(K_{\m,\Gamma})$, $y_1,\ldots y_n\in R$, and  $\a_1, \ldots, \a_n$ are integral  ideals in $\I_\m$ such that $y_i+\a_i\subsetneq \a_{\tilde{\k}}$. As explained in the proof of \cite[Proposition~6.5]{Bru1}, we can find an integral ideal $\b\in\k$ such that $\b\subsetneq \a_{\tilde{\k}}$ and $y_i+\a_i\not\subseteq \b$. The function on $\Omega$ corresponding to the projection $\prod_{i=1}^n (e_{\a_{\tilde{\k}}}-e_{y_i+\a_i})$ has $[0,\b]$ in its support, so a computation using the right-hand-side of \eqref{eqn:measures} shows that $\varphi\big(\prod_{i=1}^n (e_{\a_{\tilde{\k}}}-e_{y_i+\a_i})\big)\neq0$. 
Thus  $\pi_\varphi\big(\prod_{i=1}^n (e_{\a_{\tilde{\k}}}-e_{y_i+\a_i})\big)\neq 0$ by \cite[Corollary~8.14.4]{Ped},  so that $\pi_\varphi$ satisfies the faithfulness criteria given in \cite[Theorem~6.1]{Bru1}.\end{proof}

\begin{corollary}
\label{cor:infinitelymany}
If $K$ is not equal to $\QQ$ and is not an imaginary quadratic field, then the C*-dynamical system $(C_\lambda^*(R\rtimes R_{\m,\Gamma}),\sigma)$ has infinitely many distinct partition functions.
\end{corollary}
\begin{proof}
In light of Theorem~\ref{thm:partitionfncI}, it suffices to show that the set 
\[
\{|O| : O\subseteq \widehat{R} \text{ is a finite orbit for } R_{\m,\Gamma}^*\curvearrowright \widehat{R}\}
\]
is infinite. 

Observe that our assumption about $K$ implies that $R^*$ is infinite by Dirichlet's unit theorem. The subgroup $R_{\m,\Gamma}^*\subseteq R^*$ is of finite index, so $R_{\m,\Gamma}^*$ is also infinite. Hence, there exists an element in $ R_{\m,\Gamma}^*$ of infinite order. Using such an element and considering points of the form $\frac{1}{N}(1,1,...,1)+\ZZ^n\in\RR^n/\ZZ^n\cong\widehat{R}$ for $N\in\NN^\times$, it is not difficult to  come up with arbitrarily large orbits for $R_{\m,\Gamma}^*\curvearrowright \widehat{R}$, which shows that the above set is infinite.
\end{proof}

\section{Topological structure of extremal KMS states of type I}
\label{sec:topstructureKMS}

Given a C*-algebra $A$, we let $\ex_{I_n}T(A)$ denote the subset of $\ex T(A)$ consisting of those extremal tracial states on $A$ that are of type I$_n$. Also let $\ex_IT(A):=\bigcup_{n=1}^\infty\ex_{I_n}T(A)$. Thus, $\ex_{I_1}T(A)$ denotes the set of characters of $A$, these characters, or complex homomorphisms, are obviously pure states. We will need the following strengthening of \cite[Corollary~2.4]{Nesh} for the special case of characters in the crossed product of an action with finitely many fixed points.

\begin{lemma}
\label{lem:typeItraces&fixedpoints}
Let $G$ be a countable abelian group acting on a second countable locally compact Hausdorff space $X$. Assume that the set $\F:=\{x\in X : gx=x\text{ for all } g\in G\}$ is finite. For each $x\in \F$, there is a continuous injective map
\[
\rho_x:\widehat{G}\to \ex_{I_1}T(C_0(X)\rtimes G)
\]
such that $\rho_x(\charctr)(fu_g)=\charctr(g)f(x)$ for all $f\in C_0(X)$ and $g\in G$ where $u_g$ denotes the canonical unitary  corresponding to $g$ in the multiplier algebra of $C_0(X)\rtimes G$.
Moreover,  every character of $C_0(X)\rtimes G$ arises this way and the map
\[
\bigsqcup_{x\in\F}\rho_x:\bigsqcup_{x\in\F} \widehat{G}\to \ex_{I_1}T(C_0(X)\rtimes G)
\]
is a homeomorphism.
\end{lemma}
\begin{proof}
It follows from \cite[Corollary~2.4]{Nesh} that each $\rho_x$ is a well-defined injection of $\widehat{G}$ into $\ex T(C_0(X)\rtimes G)$. 
By \cite[Remark~2.5]{Nesh}, we see the GNS representation $\pi_{\rho_x(\charctr)}$ of $\rho_x(\charctr)$ satisfies 
\[
\pi_{\rho_x(\charctr)}(C_0(X)\rtimes G)''\cong L^2(X,\delta_x)\rtimes (G/G_x)\cong \CC
\]
for every $\charctr\in\widehat{G}$, so that the range of $\rho_x$ is contained in $\ex_{I_1}T(C_0(X)\rtimes G)$ (also, a direct calculation shows that $\rho_x(\charctr)$ is always a character).  Moreover, \cite[Remark~2.5]{Nesh} shows that $\tau\in\E T(C_0(X)\rtimes G)$ is of type I$_1$ if and only if $\tau$ corresponds via the parameterization given by \cite[Corollary~2.4]{Nesh} to a fixed point $x\in\F$ and a character $\charctr\in\widehat{G}$, that is, if and only if $\tau=\rho_x(\charctr)$. Thus, $\sqcup_{x\in\F}\rho_x$ is a bijection onto $\ex_{I_1}T(C_0(X)\rtimes G)$.

We need to show that $\rho_x$ is continuous. Suppose $(\charctr_i)_i$ is a net in $\widehat{G}$ converging to some $\charctr\in\widehat{G}$. Let $N\in\NN^\times$, $f_k\in C_0(X)$, $g_k\in G$ for $k=1,...,N$. If $\varepsilon>0$, then taking $j$ such that $i\geq j$ implies $|\charctr_i(g_k)-\charctr(g_k)|<\frac{\varepsilon}{N\max_{1\leq k\leq N} \|f_k\|_\infty}$ for all $k=1,...,N$ yields
\[
\left|\rho_x(\charctr_i)\left(\sum_{k=1}^Nf_ku_{g_k}\right)-\rho_x(\charctr)\left(\sum_{k=1}^Nf_ku_{g_k}\right)\right|\leq \sum_{k=1}^N|\charctr_i(g_k)-\charctr(g_k)| \|f_k\|_\infty<\varepsilon
\]
for all $j\geq i$. Hence, $\rho_x(\charctr_i)(\sum_{k=1}^Nf_ku_{g_k})\to \rho_x(\charctr)(\sum_{k=1}^Nf_ku_{g_k})$. Now an $\varepsilon/3$-argument shows that $\rho_x(\charctr_i)(a)\to \rho_x(\charctr)(a)$ for every $a\in C_0(X)\rtimes G$.

To show that $\sqcup_{x\in\F}\rho_x$ is a homeomorphism, it suffices to show that $\rho_{x_0}(\widehat{G})$ is open in $\ex_{I_1}T(C_0(X)\rtimes G)$ for each fixed $x_0\in\F$. Since $\widehat{G}$ is compact and each $\rho_x$ is continuous, so that $\rho_x(\widehat{G})$ is compact and hence closed. Thus, the set
\[
\bigsqcup_{x\in\F\setminus\{x_0\}}\rho_x(\widehat{G})
\]
is closed because $\F$ is assumed to be finite. Hence, its complement $\rho_{x_0}(\widehat{G})$ is open.
\end{proof}

The actions of the group of units on the duals of the additive groups of ideals representing different ideal classes need not be conjugate. Nevertheless, we have the following solidarity result for the number of fixed points of those actions. 
\begin{proposition}
\label{prop:orbits}
For each integral ideal $\a\in \I_\m$, let 
\[
\F_\a:=\{\chi \in \widehat{\a} : u\chi=\chi \text{ for all }u\in R_{\m,\Gamma}^*\}.
\]
If $R_{\m,\Gamma}^*\neq\{1\}$, then $\F_R$ is finite and $|\F_\a|=|\F_R|$ for every integral ideal $\a\in \I_\m$.
\end{proposition}
\begin{proof}
Let $u\in R_{\m,\Gamma}^*\setminus\{1\}$, and let $\a$ be an integral ideal in $\I_\m$.
Let $x_1,x_2,...,x_n$ be a $\ZZ$-basis for $\a$, where $n=[K:\QQ]$, and let $A$ denote the matrix for the $\ZZ$-linear map $\a\to\a$ given by $x\mapsto ux$ with respect to this basis. The elements $x_1,x_2,...,x_n$ are a $\QQ$-basis for $K$, so $A$ is also the matrix for the $\QQ$-linear map $K\to K$ given by $x\mapsto ux$. By \cite[Chapter~I,~\S~2,~Proposition~2.6]{Neu}, the eigenvalues for this transformation are precisely $w_1(u),w_2(u),...,w_n(u)$ where $w_1,w_2,...,w_n$ are the archimedean embeddings of $K$.
Hence, the transpose of $A$ has no eigenvalues equal to $1$, so $\{\chi\in\widehat{\a} : u\chi=\chi\}$ is finite by \cite[Lemma~2]{Hal}. This implies that $\F_\a$ is finite. 

Next we show the inequality $|\F_R|\leq |\F_\a|$. Let $F_u:=\{\chi\in\widehat{R} : u\chi=\chi\}$. By possibly replacing $\a$ with another ideal in the class $[\a]$, we may assume that $N(\a)$ and $|F_u|$ are relatively prime (such an ideal exists, for instance, by \cite[Chapter~VIII,~Theorem~7.2]{MilCFT}). Let $x_1,...,x_n$ be a $\ZZ$-basis for $R$ for which there exists $d_1,...,d_n\in\ZZ_{>0}$ be such that $d_1x_2,...,d_nx_n$ is a $\ZZ$-basis for $\a$ (such bases exist by, for example, by the theorem in \cite{MY}). Since $N(\a)=d_1\cdots d_n$, each $d_i$ is relatively prime to $|F_u|$.

Use the bases $x_1,...,x_n$ and $d_1x_1,...,d_nx_n$ to identify $\widehat{R}$ and $\widehat{\a}$, respectively, with $\RR^n/\ZZ^n$. Under these identification, the canonical projection map $\widehat{R}\to \widehat{\a}$ is taken to the map $D:\RR^n/\ZZ^n\to \RR^n/\ZZ^n$ given by $D((c_1,...,c_n)+\ZZ^n)=(d_1c_1,...,d_nc_n)+\ZZ^n$. 
Since $D$ is $R_{\m,\Gamma}^*$-equivariant, to show the inequality $|\F_R|\leq |\F_\a|$, it suffices to show that $D$ is injective on the (image of) $\F_R$. Since $\F_R\subseteq F_u$, it is even enough to show that $D$ is injective on $F_u$, where we view $F_u$ as a subgroup of $\RR^n/\ZZ^n$.
Since $F_u$ is a finite subgroup of $\RR^n/\ZZ^n$, we have $F_u\subseteq \tors(\RR^n/\ZZ^n)=\QQ^n/\ZZ^n$. In fact, since the order of every element of $F_u$ divides $|F_u|$, we have $F_u\subseteq \frac{1}{|F_u|}\ZZ^n/\ZZ^n$.
Hence, it suffices to show that $D$ is injective on $\frac{1}{|F_u|}\ZZ^n/\ZZ^n$. Suppose $(a_1,...,a_n)+\ZZ^n$ is in $\frac{1}{|F_u|}\ZZ^n/\ZZ^n$ with $D((a_1,...,a_n)+\ZZ^n)=(0,0,...,0)+\ZZ^n$. Then we have $d_ia_i\in\ZZ$ for $i=1,...,n$. For each $i$, there exists $a_i'\in \ZZ$ such that $a_i=\frac{a_i'}{|F_u|}$. Now we have $d_ia_i'\in |F_u|\ZZ$ for $i=1,...,n$. Since $d_i$ and $|F_u|$ are relatively prime, this forces $a_i'\in |F_u|\ZZ$ for $i=1,...,n$, which shows that $(a_1,...,a_n)+\ZZ^n=(0,0,...,0)+\ZZ^n$. Thus, $D$ is injective on $\frac{1}{|F_u|}\ZZ^n/\ZZ^n$, as desired.

Now we need to show that $|\F_\a|\leq |\F_R|$. First, notice that $|\F_R|=|\F_{(b)}|$ for every $b\in R_{\m,\Gamma}$. Thus, if we can find $b\in R_{\m,\Gamma}$ such that 
\begin{itemize}
\item $(b)\subseteq\a$, and
\item $\gcd([\a:(b)],|F_u|)=1$,
\end{itemize}
then we can run the same argument that we used to prove the inequality $|\F_R|\leq |\F_\a|$. For this, choose $\p\in[\a]^{-1}$ such that $N(\p)$ and $|F_u|$ are relatively prime (such a prime exists, for instance, by \cite[Chapter~VIII,~Theorem~7.2]{MilCFT}). Then $\p\a=(b)$ for some $b\in R_{\m,\Gamma}$, and $[\a:(b)]=\frac{N(b)}{N(\a)}=N(\p)$ is relatively prime to $|F_u|$.
\end{proof}

In order to select a distinguished subset from among all the partition functions associated to our systems we
will use their residues at the critical inverse temperature. We begin by computing the value of the residues at $\beta =1$ of the partial zeta functions associated to classes in $\I_\m/i(K_{\m,\Gamma})$.

\begin{proposition}
\label{prop:residue}
Let $r$ and $s$ denote the number of real and complex embeddings of $K$, respectively. If $\k\in\I_\m / i(K_{\m,\Gamma})$, then
\begin{equation}\label{eqn:residue}
\underset{\beta=1}{\Res}\zeta_\k(\beta)=|\bar{\Gamma}|\cdot \left(\frac{2^r(2\pi)^s\reg(K)[R^*:R_{\m,1}^*]}{w_\m N(\m_0)2^{r_0}|D_{K/\QQ}|^{1/2}}\right)
\end{equation}
where $\bar{\Gamma}=i(K_{\m,\Gamma})/i(K_{\m,1})$, and $r_0$ is the number of real places that divide $\m$, and $w_\m$ is the number of roots of unity in $R_{\m,1}^*$. In particular, $\underset{\beta=1}{\Res}\zeta_\k(\beta)$ does not depend on $\k$.
\end{proposition}
\begin{proof}
By \cite[Chapter~VI,~Corollary~2.9]{MilCFT}, the residue at $\beta=1$ of the partial zeta function from a class in $\Cl_\m$ is given by  
\[
\frac{2^r(2\pi)^s\reg(K)[R^*:R_{\m,1}^*]}{w_\m N(\m_0)2^{r_0}|D_{K/\QQ}|^{1/2}},
\]
and is independent of the specific class.
Let $\k\in\I_\m / i(K_{\m,\Gamma})$. Since $\I_\m / i(K_{\m,\Gamma})=\Cl_\m/(i(K_{\m,\Gamma})/i(K_{\m,1}))$, the class $\k$ is the disjoint union of $|\bar{\Gamma}|$ classes in $\Cl_\m$ all having the same residue at $\beta= 1$.  Additivity of residues then gives
$\underset{\beta=1}{\Res}\zeta_\k(\beta)=|\bar{\Gamma}|\cdot  \underset{\beta=1}{\Res}\zeta_{\tilde{\k}}(\beta)$ where $\tilde{\k}$ is any given class  in $\Cl_\m$ that is contained in $\k$. 
\end{proof}

The following lemma shows that the partial zeta function $\zeta_\k(s-1)$ can be recovered from the partition function $\unk^s\zeta_\k(s-1)$.
A function of the form $\sum_{n=1}^\infty a_nn^{-s}$ for some sequence of non-negative integers $(a_n)_{n\geq 1}$ is called a Dirichlet series. 
\begin{lemma}
\label{lem:dirichlet}
The number $\unk$ is the smallest positive integer $k$ such that $k^{-s}\unk^s\zeta_\k(s-1)$ is a Dirichlet series. 
In particular, we can recover $\unk$ from $\unk^s\zeta_\k(s-1)$.
\end{lemma}
\begin{proof}
First, we prove that if a function $f(s)$ has a series expansion $f(s)=\sum_{n=1}^\infty a_n e^{-s\lambda_n}$ for some strictly increasing sequence of real numbers $\lambda_1<\lambda_2<\cdots$ and some sequence of non-negative integers $a_1,a_2,...$, then such a series expansion is unique. 
Let $n_1<n_2<\cdots$ be the natural numbers for which the corresponding coefficients of $f(s)$ are non-zero. Then
$\lambda_{n_1}$ is the unique positive real number such that the limit $\lim_{s\to\infty}e^{s\lambda_{n_1}}f(s)$ exists and is non-zero, and
$a_{n_1}=\lim_{s\to\infty}e^{s\lambda_{n_1}}f(s).$
Then $\lambda_{n_2}$ is the unique positive real number such that the limit $\lim_{s\to\infty}e^{s\lambda_{n_2}}(f(s)-a_{n_1}e^{-s\lambda_{n_1}})$ exists and is non-zero, and $a_{n_2}=\lim_{s\to\infty}e^{s\lambda_{n_2}}(f(s)-a_{n_1}e^{-s\lambda_{n_1}})$. Continuing this way, one can compute the coefficients $a_{n_i}$ and the exponents $\lambda_{n_i}$ for $i\geq 1$. This shows the uniqueness of series expansions of $f(s)$. 

Let $f(s):=k^{-s}\unk^s\zeta_\k(s-1)$. Then, $f(s)$ has a series expansion of the form 
\begin{equation} \label{eqn:dirichlet}
f(s)=\sum_{\a\in\k, \a\subseteq R} N(\a) \left(\frac{kN(\a)}{\unk}\right)^{-s}. 
\end{equation}
Since  $\k$ contains (infinitely many) integral ideals whose norms are coprime to $\unk$, the right hand side of \eqref{eqn:dirichlet} is a Dirichlet series if and only if $\unk$ divides $k$. 
\end{proof}

\begin{definition}\label{def:minKMS}
For each $\beta>2$, let $\Sigma_{\beta,I}$ denote the set of extremal $\sigma$-KMS$_\beta$ states on $C_\lambda^*(R\rtimes R_{\m,\Gamma})$ that are of type I, and for each $\varphi\in\Sigma_{\beta,I}$, let $\tilde{Z}_\varphi(s):=N_\varphi^{-s}Z_\varphi(s)$ where $N_\varphi$ is the smallest positive integer such that $N_\varphi^{-s}Z_\varphi(s)$ is a Dirichlet series.

We say that $\varphi\in \Sigma_{\beta,I}$ is \emph{minimal} if
\[
\underset{\beta=2}{\Res}\tilde{Z}_\varphi(\beta)\leq \underset{\beta=2}{\Res}\tilde{Z}_{\varphi'}(\beta)
\]
for every $\varphi'\in \Sigma_{\beta,I}$. Let $\Sigma_{\beta,I}^{\min}$ denote the set of minimal states in $\Sigma_{\beta,I}$.
\end{definition}
By Theorem~\ref{thm:partitionfncI} and Proposition~\ref{prop:residue}, the function $\tilde{Z}_\varphi(\beta)$ associated to $\varphi\in\Sigma_{\beta,I}$ has a residue at $\beta=2$ which is an integer multiple of the positive real  number given in \eqref{eqn:residue}.

If $\varphi$ is the extremal KMS$_\beta$ state associated with the class $\k$ and a finite orbit $O$, so that Theorem~\ref{thm:partitionfncI} implies that $Z_\varphi(s)=|O| \, \unk^s \, \zeta_\k(s-1)$, then 
\[
\tilde{Z}_\varphi(s) = |O|\zeta_\k(s-1).
\]
By Lemma~\ref{lem:dirichlet}, the function $\tilde{Z}_\varphi(s)$ is explicitly computable in terms of the partition function $Z_\varphi(s)$, and conversely the partition function $Z_\varphi(s)$ is explicitly computable in terms of $\tilde{Z}_\varphi(s)$.

We define a zeta function associated to the modulus $\m$ by
 \[
 \displaystyle \zeta_{K,\m}(s):=\prod_{\p\nmid \m_0}(1-N(\p)^{-s})^{-1}=\sum_{\a\in \I_\m,\a\subseteq R}N(\a)^{-s},
 \] 
which only depends on the support of $\m_0$ and satisfies $\zeta_{K,\m}(s) = \displaystyle \sum_{\k\in \I_\m / i(K_{\m,\Gamma})}\zeta_\k(s)$.

\begin{theorem}\label{thm:main}
Let $K$ be a number field with ring of integers $R$, $\m$ a modulus for $K$, $\Gamma$ a subgroup of $(R/\m)^*$, and $R_{\m,\Gamma}$ the associated congruence monoid.  Let $\pi_0(\Sigma_{\beta,I}^{\min})$ denote the set of connected components of $\Sigma_{\beta,I}^{\min}$, and let $\bar{\varphi}$ denote the connected component of $\varphi\in\Sigma_{\beta,I}^{\min}$. Then
\begin{enumerate}
\item[\textup{(i)}] $\ \quad\quad\quad\quad\quad\quad \displaystyle |\pi_0(\Sigma_{\beta,I}^{\min})|=\begin{cases}
|\I_\m / i(K_{\m,\Gamma})| & \text{ if }R_{\m,\Gamma}^*=\{1\},\\
|\tors(R_{\m,\Gamma}^*)|\cdot |\I_\m / i(K_{\m,\Gamma})|\cdot |\F_R| & \text{ if }R_{\m,\Gamma}^*\neq \{1\};
\end{cases}$

\medskip
\item[\textup{(ii)}] the partition function $Z_\varphi$ depends only on the connected component of $\varphi$, and
\[
\sum_{\bar{\varphi}\in \pi_0(\Sigma_{\beta,I}^{\min})}\tilde{Z}_{\varphi}(s)=
\begin{cases}
\zeta_{K,\m}(s-1) & \text{ if } R_{\m,\Gamma}^*= \{1\},\\
|\tors(R_{\m,\Gamma}^*)|\cdot|\F_R| \cdot\zeta_{K,\m}(s-1) & \text{ if } R_{\m,\Gamma}^*\neq \{1\},
\end{cases} 
\]
where the sum is taken over any set of representatives and $\Re s > 2$;
\medskip\item[\textup{(iii)}] $\quad\quad\quad\quad\quad\quad \displaystyle
\underset{s\to\infty}{\lim}\sum_{\bar{\varphi}\in \pi_0(\Sigma_{\beta,I}^{\min})}\tilde{Z}_\varphi(s)=\begin{cases}
1 & \text{ if } R_{\m,\Gamma}^*=\{1\},\\
|\tors(R_{\m,\Gamma}^*)|\cdot|\F_R| & \text{ if } R_{\m,\Gamma}^*\neq \{1\}.
\end{cases}
$
\end{enumerate}
\end{theorem}
\begin{proof}

(i): By Corollary~\ref{cor:typeKMS}, and Lemma~\ref{lem:tracesonsum}, 
the correspondence from \cite[Theorem~3.2(iii)]{Bru2} specified in equation \eqref{eqn:KMSformula1} gives a homeomorphism 
\[
\Sigma_{\beta,I}\cong \bigsqcup_{\k\in \I_\m / i(K_{\m,\Gamma})}\ex_IT(C^*(\a_\k\rtimes R_{\m,\Gamma}^*)).
\]

By Theorem~\ref{thm:partitionfncI} and Proposition~\ref{prop:residue}, the KMS$_\beta$ state $\varphi=\varphi_{\beta,\k,\tau}$ is minimal if and only if $\mu_\tau$ is a probability measure on $\widehat{\a}_\k$ that is concentrated at a fixed point $\gamma\in\F_\k$, in which case $\chi_{\tau}$ can be any character of $R_{\m,\Gamma}^*$. Thus,
\[
\Sigma_{\beta,I}^{\min}\cong \bigsqcup_{\k\in \I_\m / i(K_{\m,\Gamma})}\ex_{I_1}T(C^*(\a_\k\rtimes R_{\m,\Gamma}^*)).
\]
If $R_{\m,\Gamma}^*=\{1\}$, so that $C^*(\a_\k\rtimes R_{\m,\Gamma}^*)\cong C(\widehat{\a}_\k)$ and $\ex_{I_1}T(C^*(\a_\k\rtimes R_{\m,\Gamma}^*))\cong \ex_{I_1}T(C(\widehat{\a}_\k))\cong \TT^n$, where $n=[K:\QQ]$, is connected, then we see that $\Sigma_{\beta,I}^{\min}$ has $|\I_\m / i(K_{\m,\Gamma})|$ connected components.

Now assume that $R_{\m,\Gamma}^*\neq \{1\}$. Then Proposition~\ref{prop:orbits} implies that  $\F_{\a_\k}$ is finite for each $\k$, so we can use Lemma~\ref{lem:typeItraces&fixedpoints} to get the following decomposition 
\[
\Sigma_{\beta,I}^{\min}\cong \bigsqcup_{\k\in \I_\m / i(K_{\m,\Gamma})}\ex_{I_1}T(C^*(\a_\k\rtimes R_{\m,\Gamma}^*))\cong \bigsqcup_{\k\in \I_\m / i(K_{\m,\Gamma})}\bigsqcup_{\gamma\in\F_{\a_\k}}\widehat{R_{\m,\Gamma}^*}.
\]

Since $R_{\m,\Gamma}^*$ is of finite index in $R^*$, it follows from Dirichlet's unit theorem that we have an isomorphism $R_{\m,\Gamma}^*\cong \tors(R_{\m,\Gamma}^*)\times \ZZ^d$ where $d=\rk_\ZZ(R^*)$. Hence $\widehat{R_{\m,\Gamma}^*}\cong \widehat{\tors(R_{\m,\Gamma}^*)}\times\TT^d$, so that $\widehat{R_{\m,\Gamma}^*}$ has $|\tors(R_{\m,\Gamma}^*)|$ connected components. By Proposition~\ref{prop:orbits}, $|\F_{\a_\k}|=|\F_R|$ for every $\k$, so $\Sigma_{\beta,I}^{\min}$ has 
\[
|\tors(R_{\m,\Gamma}^*)|\cdot \sum_{\k\in\I_\m/i(K_{\m,\Gamma})}|\F_{\a_\k}|=|\tors(R_{\m,\Gamma}^*)|\cdot |\I_\m / i(K_{\m,\Gamma})|\cdot |\F_R|
\]
connected components.

(ii): The partition function of $\varphi\in \Sigma_{\beta,I}^{\min}$ depends only on the class $\k$ over which $\varphi$ lies and is equal to $N(\a_\k)^s\zeta_\k(s-1)$ by Theorem~\ref{thm:partitionfncI} and Proposition~\ref{prop:orbits}.

(iii): This follows from (ii) combined with the observation that  $\lim_{s\to\infty}\zeta_{K,\m}(s-1)=1$.
\end{proof}

In the next corollary we obtain explicitly computable C*-dynamical invariants for our systems. To set the notation, 
let $\P_\QQ$ be the set of all prime numbers, and recall that the zeta function of the trivial class is $\zeta_{[R]}(s) := \sum_{(a)\in [R], a\in R_{\m,\Gamma}} N(a)^{-s}$.

\begin{corollary}\label{cor:invariants}
The following are explicitly computable invariants in terms of the partition functions of the  C*-dynamical system $(C_\lambda^*(R\rtimes R_{\m,\Gamma}),\sigma)$:
\begin{itemize}
\item[(i)] the generalized class number $|\I_\m / i(K_{\m,\Gamma})|$;
\item[(ii)] the zeta function  $\zeta_{K,\m}(s)$ associated to the modulus $\m$;
\item[(iii)] the product  $|\tors(R_{\m,\Gamma}^*)|\cdot |\F_R|$; 
\item[(iv)] the partial zeta function of the trivial class  $\zeta_{[R]}(s)$; and
\item[(v)] the subset of rational primes $\{p\in \P_\QQ: N(a) = p \text{ for some prime element } a \in  R_{\m,\Gamma}\}$.
\end{itemize}
\end{corollary}
\begin{proof}
From our analysis of type I factor KMS states and the definition of the associated partition functions we know that the family of partition functions, the space $\Sigma_{\beta,I}$ and its 
subspace $\Sigma_{\beta,I}^{\min}$ of minimal states from Definition~\ref{def:minKMS} are 
determined by the dynamical system
$(C_\lambda^*(R\rtimes R_{\m,\Gamma}), \sigma)$. Thus, the quantities in parts (i)--(iii)  of \thmref{thm:main} are preserved under equivariant isomorphism. In order to prove that the quantities listed in parts (i)--(iv) above are preserved under equivariant isomorphism, we show how to extract them from the family of zeta functions and the right hand sides of the equations in parts (i)--(iii) of \thmref{thm:main}.

By parts (i) and (iii) of \thmref{thm:main},  $|\I_\m / i(K_{\m,\Gamma})|$ is the number of connected components of $\Sigma_{\beta,I}^{\min}$ divided by the quantity $\underset{s\to\infty}{\lim}\sum_{\bar{\varphi}\in \pi_0(\Sigma_{\beta,I}^{\min})}\tilde{Z}_\varphi(s)$. This proves that $|\I_\m / i(K_{\m,\Gamma})|$ is preserved under equivariant isomorphisms.

By parts (ii) and (iii) of ~\thmref{thm:main}, we can read off $\zeta_{K,\m}(s)$ as the function $\sum_{\bar{\varphi}\in \pi_0(\Sigma_{\beta,I}^{\min})}\tilde{Z}_\varphi(s +1)$ divided by the quantity $\underset{s\to\infty}{\lim}\sum_{\bar{\varphi}\in \pi_0(\Sigma_{\beta,I}^{\min})}\tilde{Z}_\varphi(s)$. This proves that also $\zeta_{K,\m}(s)$ is  a C*-dynamical invariant.

It follows from part (iii) of \thmref{thm:main} that if $\underset{s\to\infty}{\lim}\sum_{\bar{\varphi}\in \pi_0(\Sigma_{\beta,I}^{\min})}\tilde{Z}_\varphi(s) = 1$, then $R_{\m,\Gamma}^*=\{1\}$ so that every point is a fixed point and $|\F_R| = \infty$. Thus the quantity
\[
|\tors(R_{\m,\Gamma}^*)|\cdot|\F_R|=
\begin{cases}
\underset{s\to\infty}{\lim}\sum_{\bar{\varphi}\in \pi_0(\Sigma_{\beta,I}^{\min})}\tilde{Z}_\varphi(s) & \text{ if } \underset{s\to\infty}{\lim}\sum_{\bar{\varphi}\in \pi_0(\Sigma_{\beta,I}^{\min})}\tilde{Z}_\varphi(s)>1,\\
\infty &\text{ if }\underset{s\to\infty}{\lim}\sum_{\bar{\varphi}\in \pi_0(\Sigma_{\beta,I}^{\min})}\tilde{Z}_\varphi(s) =1,
\end{cases}
\] is invariant under equivariant isomorphism. This proves that $|\tors(R_{\m,\Gamma}^*)|\cdot|\F_R|$, from part (iii) above is preserved under equivariant isomorphisms.

Suppose now $\varphi$ is a type I factor KMS$_\beta$ state. From the affine isomorphism given in \cite[Theorem~3.2(iii)]{Bru2} combined with \lemref{lem:tracesonsum}
and with the parametrization of tracial states via ergodic invariant measures in \thmref{thm:tracesandmeasures}(iii) we know that $\varphi$ determines an ideal class $\k$, and a finite orbit $O$ for the action $R_{\m,\Gamma}^*\curvearrowright\widehat{\a_\k}$. 
We have that the limit $\lim_{s \to \infty}\tilde{Z}_\varphi(s)$ is nonzero if and only if $\k$ 
is the trivial class $[R]$, in which case the limit is $|O|$.
Hence, $\zeta_{[R]}(s) = \tilde{Z}_\varphi(s+1)$ exactly for those states $\varphi$ whose partition functions satisfy $\lim_{s \to \infty}\tilde{Z}_\varphi(s) = 1$. This shows that $\zeta_{[R]}(s)$ is an invariant.

Next write $\zeta_{[R]}(s)$ as a Dirichlet series $\sum_{n=1}^\infty a_n n^{-s}$ with 
$a_n:=|\{(a) : a\in R_{\m,\Gamma}, N(a)=n\}|$, then \[
\{p\in \P_\QQ: N(a) = p \text{ for some  prime element } a \in  R_{\m,\Gamma}\} 
= \{p \in\P_\QQ: a_n \neq 0  \}
\]  is an invariant by  part (iv). This proves (v) and concludes the proof the corollary.
\end{proof}

\section{Recovering class field theoretic information}
\label{sec:recon}
The purpose of this section is to  show that the C*-dynamical invariants derived from the partition functions of type I factor KMS states of the system $(C_\lambda^*(R\rtimes R_{\m,\Gamma}),\sigma)$ produce invariants for equivariant isomorphism of systems that have class field theory content. The invariants themselves have already been obtained via K-theory, \cite[Theorem 5.5]{BruLi}, so they are invariant under  C*-algebra isomorphism without reference to dynamics. This is stronger than what we prove here using partition functions because we only show that the quantities involved are invariant under $\RR$-equivariant isomorphism. Nevertheless, we wish to reinforce the theme that KMS states deliver similar invariants to those obtained via K-theory, often through entirely different considerations. There is also a small exception in that our result on the generalized class numbers $|\I_\m/i(K_{\m,\Gamma})|$ is slightly stronger than the analogous conclusion in \cite[Theorem 5.5]{BruLi} because we do not need to assume that the class fields $K(\m)^{\bar{\Gamma}}$ and $L(\n)^{\bar{\Lambda}}$ are Galois over $\QQ$.  
We begin by recalling some number-theoretic background leading to two notions of equivalence for number fields.

Let $K$ be a number field with ring of integers $R$; in general we will let $\P_K$ denote the set of non-zero prime ideals in $R$; although in the rational case we invariably think of $\P_\QQ$ as being (represented by) the set of positive prime numbers. If $p\in\P_\QQ$ is a prime number, then $pR=\prod_{i=1}^r\p_i^{e_i}$ where the $\p_i$ are distinct prime ideals of $R$ and the $e_i$ are positive integers. The \emph{inertia degree of $\p_i$ over $p$}, denoted by $f(\p_i\vert p)$, is the degree of $R/\p_i$ as a field extension of $\ZZ/p\ZZ$ (see \cite[Chapter~I,~\S~8]{Neu}). 

The $r$-tuple $A_{K/\QQ}(p):=(f(\p_1\vert p),f(\p_2\vert p),...,f(\p_r\vert p))$, re-indexed such that $f(\p_1\vert p)\leq f(\p_2\vert p)\leq\cdots\leq f(\p_r\vert p)$, is called the \emph{splitting (or decomposition) type} of $p$ in $K$.
Let $L$ be another number field. By \cite[Theorem~1]{Per}, the following statements are equivalent:
\begin{itemize}
\item[(a)] $\zeta_K(s)=\zeta_L(s)$;
\item[(b)] $A_{K/\QQ}(p)=A_{L/\QQ}(p)$ for all $p\in\P_\QQ$;
\item[(c)] $A_{K/\QQ}(p)=A_{L/\QQ}(p)$ for all but finitely many $p\in\P_\QQ$.
\end{itemize}
Number fields $K$ and $L$ satisfying the above conditions are said to be \emph{arithmetically equivalent (over $\QQ$)}.

Our partition functions are associated to congruence monoids, and thus depend on the choice of modulus, so in order to show that number fields with the same dynamical invariants are arithmetically equivalent we need to rely on the fact that the zeta function of the modulus characterizes the zeta function of the underlying field itself; this is probably known to experts but we include a statement and proof here because we have not been able to find a suitable reference.
\begin{proposition}
\label{prop:zetam}
Let $K$ and $L$ be number fields, and let $\m$ and $\n$ be moduli for $K$ and $L$, respectively. If $\zeta_{K,\m}(s)=\zeta_{L,\n}(s)$, then $\zeta_K(s)=\zeta_L(s)$.
\end{proposition}

\begin{proof}
Suppose that $\zeta_{K,\m}(s)=\zeta_{L,\n}(s)$. By \cite[Theorem~1,~(c)$\implies$(a)]{Per}, in order to show that $\zeta_K(s)=\zeta_L(s)$, it suffices to show that $A_{K/\QQ}(p)=A_{L/\QQ}(p)$ for all but finitely many rational primes $p$. We will show that, for all but finitely many $p$, we can read off $A_{K/\QQ}(p)$ from $\zeta_{K,\m}(s)$.

Write $\zeta_{K,\m}(s)=\sum_{n=1}^\infty a_nn^{-s}$ where $a_n$ is the number of ideals of $R$ that are relatively prime to $\m_0$ and have norm $n$. The function $\zeta_{K,\m}(s)$ determines (and is determined by) the coefficients $a_n$ for $n\geq 1$. Define the support of the modulus $\m_0$ in $\QQ$ to be the finite set
\[
\supp_\QQ(\m_0):=\{p\in\P_\QQ : p\ZZ=\p\cap\ZZ\text{ for some }\p\mid\m_0\} =\{p\in\P_\QQ : \exists \p\in\P_K \text{ with } \p\mid\m_0 \text{ and } \p\mid p\}.
\] 
For  $p\in\P_\QQ\setminus \supp_\QQ(\m_0)$ and $\p\in\P_K$, we have $\p\mid p\implies \p\nmid\m_0$. Hence, for $p\in\P_\QQ\setminus \supp_\QQ(\m_0)$ and integer $f\geq 1$, we have
\[
 a_{p^f}= |\{ \a \in \I_\m:  \a\subseteq R, N(\a)=p^f\}|,
\] 
We will be done once we show that, for each $p\notin \supp_\QQ(\m_0)$, the splitting type of $p$ is determined by the numbers $a_{p^f}$ for $f\geq 1$. The idea needed to show this appears in the proof of \cite[Theorem~1]{Per}. For each prime $p$ and $f\geq 1$, let $b_{p^f}:=|\{\p\in\P_K : N(\p)=p^f\}|$. It is not difficult to see that the numbers $b_{p^f}$ for $f\geq 1$ determine the splitting type of $p$ in $K$, so we only need to show that each $b_{p^f}$ is determined from the numbers $a_{p^{f'}}$ for $f'\geq 1$. This is proven in \cite[Chapter~I,~Theorem~2.1]{Kli}. The argument  given in the proof of \cite[Theorem~1,~(a)$\implies$(b)]{Per}, also implies this, but we believe that the explicit formula expressing each $b_{p^f}$ in terms of the numbers $a_{p^{f'}}$ for $f'\geq 1$ given there is not always valid.
\end{proof}

If $\Kz$ is a number field, then the \emph{Kronecker set of $\Kz$ over $\QQ$} is
\[
D(\Kz\vert\QQ):=\{p\in\P_\QQ: \exists\mathfrak{P}\in\P_\Kz \text{ such that } \mathfrak{P}\mid p\text{ and } f_{\Kz/\QQ}(\mathfrak{P}\vert p)=1\}
\]
where $f_{\Kz/\QQ}(\mathfrak{P}\vert p)$ is the inertia degree of $\mathfrak{P}$ over $p$. Two number fields are said to be \emph{Kronecker equivalent (over $\QQ$)} if their Kronecker sets differ by only finitely many primes (see \cite{Jeh}). 
If $\Kz$ is Galois over $\QQ$, then $D(\Kz\vert\QQ)$ is, up to finitely many exceptions, the set of rational primes that split completely in $\Kz$. Since Galois extensions are determined by the primes that split in them, $D(\Kz\vert\QQ)$ determines $\Kz$ up to isomorphism when $\Kz$ is Galois over $\QQ$ (see, for instance, \cite[Chapter~V,~Theorem~3.25]{MilCFT}).

The pair $(\m,\Gamma)$ canonically gives rise to a class field (i.e., a finite abelian extension) of $K$. This is explained in detail in \cite[\S~2.3]{BruLi}, so our treatment here will be brief. For any modulus $\m$, let $K(\m)$ denote the ray class field associated with $\m$, see, e.g. \cite[Ch. V]{MilCFT}. Thus, for instance,  the ray class field $K(1)$ associated with the trivial modulus, is the Hilbert class field of $K$, and the ray class field $K(\infty)$ associated with the modulus $\m=\m_\infty=\infty$ consisting of all real embeddings of $K$ is the narrow Hilbert class field of $K$. Class field theory gives us an isomorphism $(R/\m)^*/[R^*]_\m\cong \Gal(K(\m)/K(1))$ where $[R^*]_\m$ is the image of the unit group $R^*$ in $(R/\m)^*$. Let $\bar{\Gamma}$ be the image of $\Gamma$ under the composition 
\[
(R/\m)^*\to (R/\m)^*/[R^*]_\m\cong \Gal(K(\m)/K(1)),
\]
and let $K(\m)^{\bar{\Gamma}}\subseteq K(\m)$ be the subfield of elements that are fixed by every element in $\bar{\Gamma}$. Note that $K(\m)^{\bar{\Gamma}}$ always contains $K(1)$ and $\Gal(K(\m)^{\bar{\Gamma}}/K) \cong \I_\m / i(K_{\m,\Gamma}) $. 

\begin{theorem}{\em (cf. \cite[Theorem 5.5]{BruLi})}
\label{thm:recon}
Let $K$ and $L$ be number fields with rings of integers $R$ and $S$, respectively. Let $\m$ and $\n$ be moduli for $K$ and $L$, respectively, and let $\Gamma$ and $\Lambda$ be subgroups of $(R/\m)^*$ and $(S/\n)^*$, respectively. Suppose there is an isomorphism of C*-dynamical systems 
\[
(C_\lambda^*(R\rtimes R_{\m,\Gamma}),\sigma)\cong (C_\lambda^*(S\rtimes S_{\n,\Lambda}),\sigma).
\]
Then \begin{enumerate}
\item[\textup{(i)}] $K$ and $L$ are arithmetically equivalent  
and the class fields $K(\m)^{\bar{\Gamma}}$ and $L(\n)^{\bar{\Lambda}}$ are Kronecker equivalent;
\item[\textup{(ii)}] the generalized class numbers are the same, i.e.  $|\I^K_\m / i(K_{\m,\Gamma})|=|\I^L_\n/i(L_{\n,\Lambda})|$, and $K(\m)^{\bar{\Gamma}}$ and $L(\n)^{\bar{\Lambda}}$  have the same degree over $\QQ$;
\item[\textup{(iii)}] if $K$ or $L$ is Galois, then $K=L$ in any algebraically closed field that contains both
$K$ and $L$.
\end{enumerate}
\end{theorem}
\begin{proof}
The claim about arithmetic equivalence follows from \corref{cor:invariants}(ii) and  \proref{prop:zetam}. 
 In order to prove the claim about Kronecker equivalence
it suffices to show that 
up to a finite set of primes, the Kronecker set
$D(K(\m)^{\bar{\Gamma}}\vert\QQ)$
is an invariant of our system. 

As before, let  $\supp_\QQ(\m_0):=\{p\in\P_\QQ : p\ZZ=\p\cap\ZZ\text{ for some }\p\mid\m_0\}$ and  
for each prime  $\p \nmid \m_0$, let $f(\p)$ denote the order of the class of $\p$ in $\I_\m / i(K_{\m,\Gamma})$. From 
 the proof of \cite[Theorem~5.5(i)]{BruLi} we know that 
\begin{equation}\label{eqn:Kroneckerset}
D(K(\m)^{\bar{\Gamma}}\vert\QQ)\setminus\supp_\QQ(\m_0)=(\P_\QQ\cap \{N(\p)^{f(\p)} : \p\nmid\m_0 \}) \setminus\supp_\QQ(\m_0).
\end{equation}

Write $\zeta_{[R]}(s)=\sum_{n=1}^\infty a_n n^{-s}$ with $a_n=|\{(a) : a\in R_{\m,\Gamma}, N(a)=n\}|$
as in the proof of   \corref{cor:invariants}(v). Hence to show that, up to finite sets of primes, $D(K(\m)^{\bar{\Gamma}}\vert\QQ)$ is an invariant of our system, it is enough to show that for each $p\in\P_\QQ\setminus\supp_\QQ(\m_0)$,
\[
a_p\neq 0 \quad \iff \quad p\in (\P_\QQ\cap \{N(\p)^{f(\p)} : \p\nmid\m_0 \}) \setminus\supp_\QQ(\m_0).\] 
``$\Longrightarrow$'': Suppose $p\in\P_\QQ\setminus\supp_\QQ(\m_0)$ with $a_p\neq 0$. Then there exists $a\in R_{\m,\Gamma}$ such that $N(a)=p$. Thus, $(a)=\p$ is a prime ideal such that $f_{K/\QQ}(\p\vert p)=1$ and $f(\p)=1$. Hence, $p$ lies in $\P_\QQ\cap \{N(\p)^{f(\p)} : \p\nmid\m_0 \}\setminus\supp_\QQ(\m_0)$.

``$\Longleftarrow$'': Suppose $p\in (\P_\QQ\cap \{N(\p)^{f(\p)} : \p\nmid\m_0 \}) \setminus\supp_\QQ(\m_0)
$. Then there exists a prime $\p\nmid\m_0 $ such that $p=N(\p)^{f(\p)}$. Hence, we must have $f_{K/\QQ}(\p\vert p)=1$ and $f(\p)=1$, so that $\p=(a)$ for some $a\in R_{\m,\Gamma}$. This shows that $a_p\neq 0$.

By  \corref{cor:invariants}(i), the generalized class number is an
explicitly computable invariant associated to the C*-dynamical system.
For the second assertion, recall from equation (3) of \cite[Section 2.3]{BruLi} that the degree $[K(\m)^{\bar{\Gamma}}:K]$ is given by the generalized class number $|\I^K_\m / i(K_{\m,\Gamma})|$. Since $K$ and $L$ are arithmetically equivalent by part (i) and since  arithmetically equivalent fields have the same degree over $\QQ$ by \cite[Theorem 1]{Per}, 
this finishes the proof of part (ii).

The assertion in part (iii) follows, just as in the proof of \cite[Theorem~5.5]{BruLi}(iii), from
the fact that arithmetic equivalence implies that the fields have the same Galois closure and the same degree
\cite[Theorem~1]{Per}. 
\end{proof}

\begin{remark}
We would like to emphasize that, under the stronger overall assumption of $\RR$-equivariant isomorphism between systems that is necessary for our approach,  \thmref{thm:recon} obtains similar conclusions to those in \cite[Theorem~5.5]{BruLi}.  Statement (i) is the same here as there. Our statements in parts (ii) and (iii) are slightly stronger than the corresponding ones in  \cite[Theorem~5.5]{BruLi} because they do not assume the class fields $K(\m)^{\bar{\Gamma}}$ and $L(\n)^{\bar{\Lambda}}$ are Galois over $\QQ$. Thus, the statement of  
 \cite[Theorem~5.5]{BruLi}(iv) also follows if the dynamical systems are isomorphic, because its proof only relies on parts (i)--(iii), for which we have stronger or equal versions here.

A sufficient condition for $K(\m)^{\bar{\Gamma}}$ and $L(\n)^{\bar{\Lambda}}$ to be Galois over $\QQ$ is given in \cite[Corollary~2.8]{BruLi}.
\end{remark}

For the case $K=\QQ$, we can say  more.

\begin{corollary}
Let $\m$ and $\n$ be moduli for $\QQ$ and let $\Gamma$ and $\Lambda$ be subgroups of $(\ZZ/\m)^*$ and $(\ZZ/\n)^*$, respectively. If there is an isomorphism of C*-dynamical systems 
\[
(C_\lambda^*(\ZZ\rtimes\ZZ_{\m,\Gamma}),\sigma)\cong (C_\lambda^*(\ZZ\rtimes\ZZ_{\n,\Lambda}),\sigma),
\] 
then $\langle\pm 1\rangle\ZZ_{\m,\Gamma}=\langle\pm 1\rangle\ZZ_{\n,\Lambda}$.
\end{corollary}
\begin{proof}
By Theorem~\ref{thm:recon}(i), $\QQ(\m)^{\bar{\Gamma}}$ and $\QQ(\n)^{\bar{\Lambda}}$ are Kronecker equivalent. Since $\QQ(\m)^{\bar{\Gamma}}$ and $\QQ(\n)^{\bar{\Lambda}}$ are Galois over $\QQ$, it follows from \cite[Chapter~V,~Theorem~3.25]{MilCFT} that $\QQ(\m)^{\bar{\Gamma}} =\QQ(\n)^{\bar{\Lambda}}$. Thus by \cite[Proposition~2.2]{BruLi}, we have
\begin{equation}\label{eqn:K=Q}
\langle\pm 1\rangle\cdot(\ZZ_\n\cap \ZZ_{\m,\Gamma})=\langle\pm 1\rangle\cdot(\ZZ_\m\cap \ZZ_{\n,\Lambda}).
\end{equation}
We will be done once we show that $\supp(\m_0)=\supp(\n_0)$. By Corollary~\ref{cor:invariants}(ii), $\zeta_{\QQ,\m}(s)$ is an invariant of the system $(C_\lambda^*(\ZZ\rtimes\ZZ_{\m,\Gamma}),\sigma)$. Write $\zeta_{\QQ,\m}(s)=\sum_{n=1}^\infty a_nn^{-s}$ where $a_n=|\{\a\in \I_\m :\a\subseteq R, N(\a)=n\}|$. Now one can read off $\supp(\m_0)$ as the set of rational primes $p$ such that there exists $k\in\ZZ_{>0}$ such that $a_{p^k}=0$.
\end{proof}

\end{document}